\documentclass{amsart}

\usepackage[utf8]{inputenc}
\usepackage{fullpage}
\usepackage{stmaryrd}
\usepackage{xcolor}
\usepackage{tikz}
\usepackage{amsmath}
\usepackage{amsthm}
\usepackage{amssymb}
\usepackage{multicol}
\usepackage{hyperref}
\usepackage{graphicx}
\usepackage{caption} 
\usepackage{subcaption} 

\numberwithin{equation}{section}
\newtheorem{theorem}{Theorem}[section]
\newtheorem{corollary}[theorem]{Corollary}
\newtheorem{example}[theorem]{Example}
\newtheorem{lemma}[theorem]{Lemma}
\newtheorem{proposition}[theorem]{Proposition}
\newtheorem{problem}{Problem}
\theoremstyle{definition}
\newtheorem{definition}[theorem]{Definition}
\newtheorem{remark}[theorem]{Remark}

\newtheorem{assumption}[theorem]{Assumption}

\newcommand{\comments}[1]{}

\title{Optimal transport paths with capacity induced cost function.}
\author{Haotian Sun}
\address{
Yau Mathematical Sciences Center\\
Tsinghua University \\
Beijing, 100084, China }
\email{hatsun@mail.tsinghua.edu.cn}
\author{Qinglan Xia}
\address{
Department of Mathematics\\
University of California
at Davis\\
Davis, CA, 95616, USA }
\email{qlxia@ucdavis.edu}
\date{}

\begin{document}

\begin{abstract}
This article generalizes the study of ramified optimal transport with capacity constraint in \cite{capacity_components} by generalizing the $\mathbf{M}_{\alpha}$ cost to $\mathbf{M}_{\alpha,c}$, which incorporates capacity constraints into the cost function.
Equipped with $\mathbf{M}_{\alpha,c}$ cost, we prove the existence of optimal transport path, $\mathbf{M}_{\alpha,c}$ related inequalities, decomposition of any general transport paths, and occurrence of direct line segments in an optimal transport path.

\end{abstract}

\maketitle

\section{Introduction}
Optimal transport problem aims at finding a cost efficient way to transport mass from sources to targets, where the sources and targets are often characterized using measures.
The Monge-Kantorovich transport problem \cite{Luigi, villani2} uses transport map and transport plan to characterize the transportation between measures. 
The total cost in the Monge-Kantorovich transport problem is formulated using sources and targets, which means that it is independent of the actual ``path'' that connects any source point to target point.

Unlike the Monge-Kantorovich transport problem, the ramified (or branched) optimal transport problem \cite{xia2003,xia2015} uses transport path to characterize transportation. 
A transport path is defined using weighted directed graphs and generalized using rectifiable $1-$currents. 
Moreover, the total cost in the ramified transport problem is dependent on the ``path'' that conducted the transportation.

From \cite{xia2003} and \cite{xia2015}, we have the following definitions of ramified transportation. 
Let $X$ be a convex compact set in $\mathbb{R}^m$, an atomic measure defined on $X$ is
$$\sum_{i=1}^k m_i \delta_{x_i},$$
with distinct points $x_i \in X$, $m_i >0$ for $i=1,2,\ldots,k$. Here, $k$ can be $+\infty$.
If we further assume $k< \infty$, we may call the above atomic measure, finitely atomic measure.

Given two atomic measures,
\begin{equation}
\label{eq: atomic measures}
\textbf{a}= \sum_{i=1}^k m_i \delta_{x_i}  \ , \ \textbf{b}= \sum_{j=1}^\ell n_j \delta_{y_j},
\end{equation}
supported on $X$ of equal total mass. 
A transport path from $\textbf{a}$ to $\textbf{b}$ is a weighted directed graph $G=[V(G), E(G), w]$ consisting of a vertex set $V(G)$, a directed edge set $E(G)$ and a weight function 
$w: E(G) \rightarrow (0, +\infty)$ 
such that $\{x_1,x_2,\cdots,x_k \}\cup\{y_1,y_2,\cdots,y_\ell \} \subseteq V(G) $ and for any vertex $v\in V(G)$, there is a balance equation: 
\begin{equation}
\label{eqn: balanced_equation}
\sum_{ \substack{e\in E(G) \\  e^-=v }} w(e) \ = 
\sum_{ \substack{e\in E(G) \\  e^+=v}} w(e) \  +\  \left\{
    \begin{array}{ll}
        \ \ m_i  &\mbox{if } v=x_i \mbox{ for some } i=1,\cdots,k  \\
        -n_j     &\mbox{if }v=y_j \mbox{ for some } j=1,\cdots,\ell \\
        \ \ 0    &\mbox{otherwise}
    \end{array}
    \right.
\end{equation}
where $e^-$ and $e^+$ denote the starting and ending point of the edge $e\in E(G)$. Note that the condition (\ref{eqn: balanced_equation}) means that masses are conserved at every interior vertex.
We denote the set of all transport paths from $\mathbf{a}$ to $\mathbf{b}$ as $$Path(\mathbf{a},\mathbf{b}).$$

For any real number $\alpha \in [0,1]$, the $\mathbf{M}_\alpha$ cost  of 
$G =[V(G),E(G), w] $ 
is defined by
\begin{equation}
\label{eq: ramified transport cost function}
\textbf{M}_\alpha(G):=\sum_{e\in E(G)} \left(w(e)\right)^\alpha \mathcal{H}^1(e),
\end{equation}
where
$\mathcal{H}^1(e)$ is the
1-dimensional Hausdorff measure
or length of the edge $e$.
We say $G^* \in Path(\mathbf{a},\mathbf{b})$ is an optimal transport path if $$\mathbf{M}_{\alpha}(G^*) \le \mathbf{M}_{\alpha}(G),  \text{ for all } G \in Path(\mathbf{a},\mathbf{b}).$$

This article studies the behavior of capacity constraint on transport paths, which also serves as a continuation of the previous work \cite{capacity_components}.
In this article we generalize the $\mathbf{M}_{\alpha}$ cost used in \cite{xia2003, capacity_components} by $\mathbf{M}_{\alpha,c}$, which extends the number of admissible transport paths as compared to the admissible transport multi-path used in \cite{capacity_components}. 
In Section 2, we review some concepts related to geometric measure theory that will be used in ramified transport problem.
In Section 3, we give the existence of optimal transport path under $\mathbf{M}_{\alpha,c}$ cost, and some inequalities related to $\mathbf{M}_{\alpha,c}$ cost.
In Section 4, we show transport paths with certain overlapping property can be decomposed into the sum of  transport paths with weight equals integer multiple of the capacity constraint and a transport path with weight less than the capacity constraint. 
Also, we showed that paths in an optimal transport path with weight equals integer multiple of the capacity constraint often transport directly via some line segments.

\section{Preliminaries}
We first recall some basic concepts from geometric measure theory \cite{Simon, Lin}.

For any open set $U$ in $\mathbb{R}^{m}$
and $k \le m$, let $\mathcal{D}^k(U)$ be the set of all $C^\infty$ $k$-forms in $U$ with compact supports. The space $\mathcal{D}_k(U)$ of \textit{$k$-currents} 
is the dual space of $\mathcal{D}^k(U)$.

For any current $T \in \mathcal{D}_k(U)$, the \textit{mass} of $T$ is defined by
$$
\mathbf{M}(T)= \sup \{ T(\omega) \,:\,  \sup_{x\in U}\|\omega(x)\| \le 1, \omega \in \mathcal{D}^k(U)\}, $$
where the {\it comass}
$\|\omega(x)\|:=  \sup\{\vert\langle \omega(x), \xi \rangle\vert: \xi \text{ is a unit, simple, k-vector in } \mathbb{R}^{m}\}.$  Also, its \textit{boundary} $\partial T \in \mathcal{D}_{k-1}(U)$ is defined by 
$$\partial T(\omega) := T(d\omega), \forall  \omega \in \mathcal{D}^{k-1}(U), \text{ when } k\ge 1 ,$$
and $\partial T :=0$ when $k=0$.

A current $T\in \mathcal{D}_k(U)$ is said to be \textit{normal} if 
$\mathbf{M}(T) + \mathbf{M}(\partial T) < \infty$.
In \cite{Paolini}, Paolini and Stepanov introduced the concept of subcurrents: For any $T, S\in \mathcal{D}_k(U)$, $S$ is called a \textit{subcurrent} of $T$ if
\begin{equation}\label{eqn: subcurrent}
 \mathbf{M}(T-S) + \mathbf{M}(S) = \mathbf{M}(T).    
\end{equation}
A normal current $T\in \mathcal{D}_k(\mathbb{R}^m)$ is \textit{acyclic} if there is no non-trivial subcurrent $S$ of $T$ such that $\partial S =0$. 
Besides this acyclic definition, there is another kind of ``acyclic'' characterization of transport path, which we will call it \textit{cycle-free}.

The above characterization of ramified transportation can be generalized to transportation between Radon measures. 
Given two Radon measures $\mu^-, \mu^+$ of equal mass, both supported on $X$. 
A transport path from $\mu^-$ to $\mu^+$ is a
rectifiable 1-current $T=\underline{\underline{\tau}}(M,\theta(x),\xi(x))$ with $\partial T = \mu^+-\mu^-$. 
We denote the set of all such transport path as 
$$Path(\mu^-,\mu^+).$$
The corresponding $\mathbf{M}_{\alpha}$  cost, for $\alpha \in [0,1]$, is 
$$\mathbf{M}_{\alpha}(T) := \int_{M} \theta(x)^\alpha \,d\mathcal{H}^1. $$
A transport path $T^* \in Path(\mu^-,\mu^+)$ is optimal if $$\mathbf{M}_{\alpha}(T^*) \le \mathbf{M}_{\alpha}(T), \text{ for all } T \in Path(\mu^-,\mu^+).$$
In other words, an $\mathbf{M}_{\alpha}$ minimizer in $Path(\mu^-,\mu^+)$ is called an optimal transport path from $\mu^-$ to $\mu^+$.

\begin{definition}(\cite[Definition 4.2]{decomposition})
\label{def: cycle-free current}
Let $T=\underline{\underline{\tau}}(M,\theta,\xi)$ and $S=\underline{\underline{\tau}}(N,\phi,\zeta)$ be two real rectifiable $k$-currents.
\begin{itemize}
    \item[(a)] We say $S$ is
{\it on} $T$ if $\mathcal{H}^k(N\setminus M)=0$, and $\phi(x)\le \theta(x)$ for $\mathcal{H}^k$ almost every $ x\in N$. 
\item[(b)] $S$ is called a {\it cycle} on $T$ if $S$ is on $T$ and $\partial S=0$.
\item[(c)] $T$ is called {\it cycle-free} if except for the zero current, there is no other cycles on $T$.
\end{itemize} 
\end{definition}

In \cite{capacity_components}, we consider the transport problem with suitable capacity constraint.
Given two Radon measures $\mu^-,\mu^+$ supported on $X$ of equal mass and a capacity constraint $c>0$, 
we want to find an optimal transport path among all $T = \underline{\underline{\tau}}(M,\theta(x),\xi(x)) \in Path(\mu^-,\mu^+)$, such that 
$$\theta(x) \le c \text{ for all } x \in M.$$
However, such ``seemingly'' natural definition of transport paths with capacity fail to converge under the condition of $\theta(x) \le c$, which is demonstrated in \cite[Example $1.2$]{capacity_components}.
In order to deal with this non-convergence issue, the transport problem with capacity constraint is reformulated using ``multi-path'' in \cite{capacity_components} as follows.

\begin{definition}
Let $\mu^-$, $\mu^+$ be two Radon measures on $ \mathbb{R}^m$ with equal mass, supported on compact sets, $\alpha \in 
[0,1]$, and $c>0$. 
Minimize  
\begin{equation}\label{eqn: cost with multi-path}
\mathbf{M}_\alpha(\vec{T}):=\sum_{k=1}^\infty \mathbf{M}_\alpha(T_k)    
\end{equation}
among all
$\vec{T}=(T_1,T_2,\cdots,T_k,\cdots)$ such that for each $k$,
\begin{equation}\label{eqn: multi-path}
T_k \in Path(\mu^-_k, \mu^+_k), \   
\sum_{k=1}^\infty \mu_k^- = \mu^-, \ \sum_{k=1}^\infty \mu_k^+ = \mu^+, and\ 0<\| \mu_k^-\| = \| \mu_k^+\| \le c.
\end{equation}
Each $\vec{T}=(T_1,T_2,\cdots,T_k,\cdots)$
satisfying (\ref{eqn: multi-path}) is called a {\it transport multi-path} from $\mu^-$ to $\mu^+$ with capacity $c$. 
The family of all such transport multi-paths is denoted by $Path_c(\mu^-,\mu^+)$.
\end{definition}

Nevertheless, there are still drawbacks when characterizing transport paths with capacity using transport multi-paths.
As illustrated in \cite[Remark $1.3$, Figure $3$]{capacity_components}, there are admissible transport paths with weight on each edge less or equal to $c$ and its boundary equals the sum of boundaries of transport multi-paths, such that its $\mathbf{M}_{\alpha}$ cost is less or equal to the $\mathbf{M}_{\alpha}$ cost of transport multi-path.

Therefore, we need to update the characterization of branched transport problem with capacity constraints such that the set of admissible transport paths is ``extended" as compared to transport multi-paths, while still in some sense ``contain'' the condition, $\theta(x) \le c$. 
The approach in this paper is generalize the regular $\mathbf{M}_\alpha$ cost to $\mathbf{M}_{\alpha,c}$ cost which is defined in next section.

In the following, we present some key results from \cite{lower semicontinuous} which characterize the lower semi-continuity of suitable cost functions used in transport path problems, and it will be used to prove the existence of optimal transport paths.

Given any $m-$dimensional current $T \in \mathcal{D}_m(\mathbb{R}^n)$, the flat norm on $\mathcal{D}_m(\mathbb{R}^n)$ is defined as:
$$\mathbb{F}(T):= \inf \{\mathbf{M}(S) + \mathbf{M}(T- \partial S) : S \in \mathcal{D}_{m+1}(\mathbb{R}^n)\} .$$

\begin{assumption}\label{assump: cost function}
Consider the Borel functions $H : \mathbb{R} \to [0,\infty)$ that satisfy the following conditions.

\begin{itemize}
    \item[(1)] $H(0)=0$ and $H$ is even, namely $H(-\theta) = H(\theta)$ for every $\theta \in \mathbb{R}$;

    \item[(2)] $H$ is subadditive, namely $H(\theta_1 + \theta_2) \le H(\theta_1) + H(\theta_2)$ for every $\theta_1, \theta_2 \in \mathbb{R}$.

    \item[(3)] $H$ is lower semicontinuous, namely $H(\theta) \le \liminf_{j \to \infty} H(\theta_j)$ whenever $\theta_j$ is a sequence of real numbers such that $|\theta - \theta_j| \searrow 0$ when $j \uparrow \infty$. 
\end{itemize}

\end{assumption}

Let $H$ be as in Assumption \ref{assump: cost function} and let $\mathcal{R}_m(\mathbb{R}^n)$ be the set of rectifiable $m-$currents.  
We define, $\mathbb{M}_H$, the $H-$mass on $\mathcal{R}_m(\mathbb{R}^n)$ as 
$$\mathbb{M}_H(R):= \int_E  H(\theta(x)) d\mathcal{H}^m(x), \ \text{ for every } R=\underline{\underline{\tau}}(E,\theta,\xi)\in \mathcal{R}_m(\mathbb{R}^n) .$$

\begin{proposition}\label{prop: lower semicontinuity of cost function}

Let $H$ satisfies Assumption \ref{assump: cost function}, and let $U \subset \mathbb{R}^n$ be open. Let $T_j,T \in \mathcal{R}_m(\mathbb{R}^n)$ be rectifiable $m-$currents such that $\mathbb{F}(T-T_j) \searrow 0$ as $j \to \infty$. Then
$$\mathbb{M}_{H}(T \lfloor U) \le \liminf_{j \to \infty} \mathbb{M}_H (T_j \lfloor U).$$
    
\end{proposition}

Moreover, from \cite{Lin, Simon}, we have the following result that connects flat norm convergence with weak convergence of rectifiable currents.

\begin{proposition}\label{prop: weak convergence = flat norm convergence}
Let $T_i, T$ be rectifiable currents with 
$$\sup \{\mathbf{M}(T_i) + \mathbf{M}(\partial T_i)\}  < \infty,$$
then $T_i \rightharpoonup T$ if and only if $\mathbb{F}(T_i -T) \to 0$.
\end{proposition}

\section{Transport problem with capacity related cost function.}

\begin{definition}
Given two Radon measures $\mu^-$, $\mu^+$ supported on compact sets in $\mathbb{R}^m$ with equal total mass, $\alpha \in [0,1]$, and $c >0$. For any $T = \underline{\underline{\tau}}(M,\theta(x),\xi(x)) \in Path(\mu^-,\mu^+)$, the transport cost of $T$ is defined as:
\begin{equation}\label{eqn: cost with capacity}
\mathbf{M}_{\alpha,c}(T):= c^\alpha\cdot \int_{M} \left\lfloor\frac{\theta(x)}{c}\right\rfloor  + \left( \frac{\theta(x)}{c} - \left\lfloor\frac{\theta(x)}{c}\right\rfloor \right)^\alpha d\mathcal{H}^1,
\end{equation}
where $\left\lfloor\theta(x)/c \right\rfloor$ denotes the largest integer less or equal to $\theta(x)/c$.
\end{definition}

Note that  $\mathbf{M}_{1,c}(T)=\mathbf{M}(T)$ and  $\lim_{c\rightarrow \infty}\mathbf{M}_{\alpha, c}(T)\rightarrow \mathbf{M}_{\alpha}(T)$.

We now consider the following Plateau-type problem:
\begin{problem}
    Minimize $\mathbf{M}_{\alpha,c}(T)$ among all  $T \in Path(\mu^-, \mu^+)$, namely among all rectifiable 1-currents $T$ with $\partial T=\mu^+-\mu^-$.
\end{problem}

Note that the cost function defined in equation (\ref{eqn: cost with capacity}) allows overlapping as compared to the definition of multi-paths components in (\ref{eqn: cost with multi-path}), while implicitly restricting the maximum weight less or equal to $c$ when calculating the total cost.
Moreover, the integer $\lfloor \theta(x) /c \rfloor$ indicates that the ``total'' weight at each point is subdivided into components such that each component has weight less or equal to $c$.

In the following, we first show some preparatory works for the existence result of optimal transport path under the cost function (\ref{eqn: cost with capacity}), and then the existence result. 

Given $\alpha \in [0,1]$, $c >0$, we first consider properties of the function
\begin{equation}\label{eqn: H cost funtion}
H_{c,\alpha}(x):= \left\lfloor \frac{x}{c} \right\rfloor + \left( \frac{x}{c} - \left\lfloor \frac{x}{c} \right\rfloor \right)^\alpha  
\end{equation} 
on $\mathbb{R}$. Clearly, it has the following properties:
\begin{itemize}
    \item $H_{c,\alpha}(x)=H_{1,\alpha}(\frac{x}{c})$, where $H_{1,\alpha}(x):= \left\lfloor x \right\rfloor + \left( x - \left\lfloor x \right\rfloor \right)^\alpha$;
    \item $H_{c,\alpha}(nc)=H_{1,\alpha}(n)=n$ for each integer $n$. In particular, $H_{c,\alpha}(0) =0$;
    \item when $\alpha\in (0,1]$, $H_{c,\alpha}(x)$ is strictly increasing, concave and  continuous on $\mathbb{R}$;
    \item when $\alpha=0$, $H_{c,\alpha}(x)$ is increasing, piecewise constant, and has jump discontinuities at integers, and lower semicontinuous;

   \item For fixed $x$ and $c$, $H_{c,\alpha}(x)$ is a decreasing function of $\alpha\in[0,1]$. In particular, 
   \begin{equation}
   \label{eqn: H_comp}
   \frac{x}{c}=H_{c,1}(x)\le H_{c,\alpha}(x) \le H_{c,0}(x) \le  \left\lfloor \frac{x}{c} \right\rfloor+1.
   \end{equation}
   
   \item Follows from concavity, $H_{c,\alpha}(x)$ is subadditive in the sense that $H_{c,\alpha}(x_1+x_2)\le H_{c,\alpha}(x_1)+H_{c,\alpha}(x_2)$.
\end{itemize}

\comments{
\begin{lemma}
Given $\alpha \in [0,1]$, $c >0$, then 
\begin{equation}\label{eqn: H cost funtion}
H_c(x):= \left\lfloor \frac{x}{c} \right\rfloor + \left( \frac{x}{c} - \left\lfloor \frac{x}{c} \right\rfloor \right)^\alpha,    
\end{equation} 
is continuous for all $x\in \mathbb{R}$, with $H_c(0) =0$. Moreover, $H_c(x)=H_1(\frac{x}{c})$.
\end{lemma}

\begin{proof}
Given $x \in \mathbb{R}$, and without loss of generality assume there exists $N \in \mathbb{Z}$ such that $$N \le \frac{x}{c} <N+1,$$
then $$H_c(x) = N + \left(\frac{x}{c} - N\right)^\alpha,$$
which is continuous when $x/c \not= N$.

Next, suppose $$\frac{x}{c} = N,$$
then $$\lim_{x/c \to N^+} H_c(x) = \lim_{x/c \to N^+} N + \left(\frac{x}{c} - N\right)^\alpha = N,$$
$$\lim_{x/c \to N^-} H_c(x) = \lim_{x/c \to N^-} N-1 + \left(\frac{x}{c} - (N-1)\right)^\alpha = N-1 + 1 =N.$$
Hence, we have $H_c(x)$ is continuous on $\mathbb{R}$.
\end{proof}
}

We now discuss some properties of the $\mathbf{M}_{\alpha,c}$ cost starting from its subadditivity.

\begin{proposition} 
Let $c>0$, $\alpha \in [0,1]$. For any two rectifiable 1-currents $T_1$ and $T_2$, we have
\begin{equation}\label{eqn: M_alpha,c subadditivity}
\mathbf{M}_{\alpha,c}(T_1+T_2) \le \mathbf{M}_{\alpha,c}(T_1) + \mathbf{M}_{\alpha,c}(T_2).   
\end{equation}

\end{proposition}

\begin{proof}
    Suppose $T_k=\underline{\underline{\tau}}(M_k,\theta_k(x),\xi_k(x))$
with $k=1,2$. Then $T_1+T_2=\underline{\underline{\tau}}(M_1\cup M_2,\theta(x),\xi(x))$ is still a rectifiable 1-current with
\begin{equation*} 
\theta(x)=\theta_1(x)\text{ and } \xi(x)=\xi_1(x), \ \text{ if }x\in M_1\setminus M_2,
\end{equation*}
\begin{equation*} 
\theta(x)=\theta_2(x) \text{ and } \xi(x)=\xi_2(x), \ \text{ if }x\in M_2\setminus M_1.
\end{equation*}
On the intersection $M_1\cap M_2$, by the uniqueness of the approximate tangent line, for $\mathcal{H}^1- a.e. \ x\in M_1\cap M_2$, it has $\xi_1(x)=\pm \xi_2(x)$. Thus, $\theta(x)=\theta_1(x)+\theta_2(x)$ if $\xi_1(x)=\xi_2(x)$ and $\theta(x)=|\theta_1(x)-\theta_2(x)|$ if $\xi_1(x)=-\xi_2(x)$. In any case, it has $\theta(x)\le \theta_1(x)+\theta_2(x)$, for $\mathcal{H}^1-a.e.\ x\in M_1\cap M_2$. As a result, by the monotonicity and subadditivity of the function $H_{c,\alpha}$, it follows that
\begin{eqnarray*}
    &&\frac{1}{c^\alpha} \cdot\mathbf{M}_{\alpha, c}(T_1+T_2)=\int_{M_1\cup M_2}H_{c,\alpha}(\theta(x))d\mathcal{H}^1(x)\\
    &\leq & \int_{M_1\setminus M_2}H_{c,\alpha}(\theta_1(x))d\mathcal{H}^1(x)+\int_{ M_2\setminus M_1}H_{c,\alpha}(\theta_2(x))d\mathcal{H}^1(x)+\int_{M_1\cap M_2}H_{c,\alpha}(\theta_1(x)+\theta_2(x))d\mathcal{H}^1(x)\\
     &\leq & \int_{M_1\setminus M_2}H_{c,\alpha}(\theta_1(x))d\mathcal{H}^1(x)+\int_{ M_2\setminus M_1}H_{c,\alpha}(\theta_2(x))d\mathcal{H}^1(x)+\int_{M_1\cap M_2}H_{c,\alpha}(\theta_1(x))+H_{c,\alpha}(\theta_2(x))d\mathcal{H}^1(x)\\
     &=& \int_{M_1}H_{c,\alpha}(\theta_1(x))d\mathcal{H}^1(x)+\int_{M_2}H_{c,\alpha}(\theta_2(x))d\mathcal{H}^1(x)\\
     &=&\frac{1}{c^\alpha} \cdot \left(\mathbf{M}_{\alpha, c}(T_1)+\mathbf{M}_{\alpha, c}(T_2)\right)
\end{eqnarray*}
as desired.
\end{proof}

\begin{lemma}\label{lemma: M_alpha,c inequality}
Given $\alpha\in [0,1]$ and $c>0$, for any rectifiable 1-current $T$,
it holds that
\begin{equation}
\label{eqn: mass_comp}
    \mathbf{M}(T)\le c^{1-\alpha}\mathbf{M}_{\alpha, c}(T)\le \mathbf{M}(T)+c \cdot\mathbf{size}(T).
\end{equation}

\begin{proof}
    Let $T =\underline{\underline{\tau}}(M,\theta(x),\xi(x)) $. By (\ref{eqn: H_comp}), we have
     \[\int_M \frac{\theta(x)}{c}d \mathcal{H}^1(x)\le \int_M H_{c,\alpha}\left(\frac{\theta(x)}{c}\right)d \mathcal{H}^1(x) \le \int_M \left(\frac{\theta(x)}{c}+1 \right)d \mathcal{H}^1(x).\]
    Using $\mathbf{size}(T):=\mathcal{H}^1(M)$, it gives
    \[    \frac{\mathbf{M}(T)}{c}\le \frac{\mathbf{M}_{\alpha, c}(T)}{c^\alpha}\le  \frac{\mathbf{M}(T)}{c}+\mathbf{size}(T),\]
   yielding (\ref{eqn: mass_comp}).
\end{proof}
\comments{
Given $T\in Path(\mu^-, \mu^+)$, $\alpha\in [0,1]$, and $c>0$, it holds that
\[ \mathbf{M}_{\alpha,c}(T)\ge \frac{\mathbf{M}(T)}{c^{1-\alpha}}.\]
}
\end{lemma}

\comments{
\begin{proof}
By definition $$  \frac{\theta(x)}{c} - \left\lfloor\frac{\theta(x)}{c}\right\rfloor <1.$$
Since $\alpha \in [0,1]$, we have $$  \frac{\theta(x)}{c} - \left\lfloor\frac{\theta(x)}{c}\right\rfloor 
\le  \left( \frac{\theta(x)}{c} - \left\lfloor\frac{\theta(x)}{c}\right\rfloor \right)^{\alpha},$$
so that
\begin{eqnarray*}
\mathbf{M}(T)
&=&
\int_{M} \theta(x) d\mathcal{H}^1(x)  \\
&=&
c \cdot \int_{M} \frac{\theta(x)}{c} d\mathcal{H}^1 \\
&=&
c \cdot \int_{M} \left\lfloor\frac{\theta(x)}{c}\right\rfloor +  \left( \frac{\theta(x)}{c} - \left\lfloor\frac{\theta(x)}{c}\right\rfloor \right) d\mathcal{H}^1 \\
&\le&
c \cdot \int_{M} \left\lfloor\frac{\theta(x)}{c}\right\rfloor +  \left( \frac{\theta(x)}{c} - \left\lfloor\frac{\theta(x)}{c}\right\rfloor \right)^\alpha d\mathcal{H}^1 \\
&=&
c^{1-\alpha}\cdot c^{\alpha} \int_{M} \left\lfloor\frac{\theta(x)}{c}\right\rfloor +  \left( \frac{\theta(x)}{c} - \left\lfloor\frac{\theta(x)}{c}\right\rfloor \right)^\alpha d\mathcal{H}^1 \\
&=&
c^{1-\alpha} \mathbf{M}_{\alpha,c}(T).
\end{eqnarray*}
\end{proof}

}

\begin{theorem}
Given two Radon measures $\mu^-$, $\mu^+$ supported on compact sets in $\mathbb{R}^m$ with equal total mass.
For any $\alpha \in [0,1]$, and $c >0$, there exists $T^* \in Path(\mu^-,\mu^+)$, such that 
$$\mathbf{M}_{\alpha,c}(T^*) \le \mathbf{M}_{\alpha,c}(T),$$
for all $T \in Path(\mu^-,\mu^+)$.
\end{theorem}

\begin{proof}
By \cite{xia2003}, we have $Path(\mu^-,\mu^+)$ is non-empty.
If $$\inf\{\mathbf{M}_{\alpha,c}(T) : T \in Path(\mu^-,\mu^+)\}=\infty,$$ then any $T$ can be a $\mathbf{M}_{\alpha,c}$ minimizer.

Now, suppose $$\inf\{\mathbf{M}_{\alpha,c}(T) : T \in Path(\mu^-,\mu^+)\} < \infty.$$ 
Let $$\{T_n = \underline{\underline{\tau}}(M_n,\theta_n(x),\xi(x)) \in Path(\mu^-,\mu^+) : n\in \mathbb{Z}^+\}$$ be an arbitrary $\mathbf{M}_{\alpha,c}$-minimizing sequence, where $\mathbf{M}_{\alpha,c}(T_n) \le K$ for each $n$ and for some $K >0$.

By Lemma \ref{lemma: M_alpha,c inequality} we have 
$$\mathbf{M}(T_n) \le  c^{1-\alpha} \mathbf{M}_{\alpha,c}(T_n) \le c^{1-\alpha} K < \infty,$$
for all $n \in \mathbb{Z}^+$.
Moreover, by definition of Radon measure, we have 
$$\mathbf{M}(\partial T_n) = \|\mu^-\|(\mathbb{R}^m) +  \|\mu^+\|(\mathbb{R}^m) < \infty.$$
Therefore, there exists $T^* = \underline{\underline{\tau}}(M,\theta(x),\xi(x))$, such that $T_n \rightharpoonup T^*$.
By Proposition \ref{prop: weak convergence = flat norm convergence}, we have $$\mathbb{F}(T_n - T^*) \to 0.$$
Since $H_{c,\alpha}(x)$ is lower semicontinuous for all $x \in \mathbb{R}$ and subadditive, we may define $\tilde{H}(x) := H_{c,\alpha}(|x|)$ for all $x \in \mathbb{R}$, so that $\tilde{H}(x)$ satisfies Assumption \ref{assump: cost function}, and $$\mathbf{M}_{\alpha,c}(T) = \int_{M} \tilde{H}(\theta(x)) d \mathcal{H}^1 = \mathbb{M}_{\tilde{H}}(T).$$
By Proposition \ref{prop: lower semicontinuity of cost function}, we have 
$$\mathbf{M}_{\alpha,c}(T^*) =\mathbb{M}_{\tilde{H}}(T^*) \le \liminf_{n \to \infty} \mathbb{M}_{\tilde{H}}(T_n) = \liminf_{n \to \infty} \mathbf{M}_{\alpha,c}(T_n),$$
which implies $T^*$ is a transport path with minimum $\mathbf{M}_{\alpha,c}$ cost.
\end{proof}

\begin{proposition} 
Suppose $T \in Path(\mu^-,\mu^+)$, where $T = \underline{\underline{\tau}}(M,\theta(x),\xi(x))$, and let $c>0$, $\alpha \in [0,1]$, then 
\begin{equation}\label{eqn: M_alpha,c scalar multiple}
\mathbf{M}_{\alpha,c}(h T) \le  h^\alpha \mathbf{M}_{\alpha,c}(T), \text{ when } 0\le h\le 1, \text{ and }  h^\alpha \mathbf{M}_{\alpha,c}(T)\le \mathbf{M}_{\alpha,c}(h T), \text{ when } h \ge 1.    
\end{equation}
Note that we adopt the notation where $hT:= \underline{\underline{\tau}}(M,h \cdot \theta(x),\xi(x))$.

\end{proposition}

\begin{proof}
By definition in (\ref{eqn: cost with capacity}), we have 
$$\mathbf{M}_{\alpha,c}(h T) = c^\alpha\int_{M} \left\lfloor \frac{h\theta(x)}{c}\right\rfloor + \left( \frac{h\theta(x)}{c}-\left\lfloor \frac{h\theta(x)}{c}\right\rfloor\right)^\alpha d\mathcal{H}^1,$$
and 
$$h^\alpha \mathbf{M}_{\alpha,c}(T)= c^\alpha h^\alpha  \int_{M} \left\lfloor \frac{\theta(x)}{c}\right\rfloor + \left( \frac{\theta(x)}{c}-\left\lfloor \frac{\theta(x)}{c}\right\rfloor\right)^\alpha d\mathcal{H}^1 .$$
Therefore, it suffices to prove the inequality between the following two functions given certain values of $h$:
$$f_h(x) := h^\alpha \lfloor x \rfloor + h^\alpha (x- \lfloor x \rfloor)^\alpha \text{ and } g_h(x) := \lfloor hx \rfloor  + (hx - \lfloor hx \rfloor)^\alpha, \text{ for } x \ge 0.$$

When $h=0$ and identifying $0^0:=1$, we have $$f_0(x)=0^\alpha \lfloor x \rfloor +  0^\alpha (x- \lfloor x \rfloor)^\alpha ,\  g_0(x)=0 +(0-0)^\alpha = 0^\alpha .$$
Since $x\ge 0$, then for $0\le \alpha \le 1$, we have $f_0(x) \ge g_0(x)$.

When $h\not= 0$ and $\alpha=1$, we have 
\begin{equation}\label{eqn: alpha=1 equality}
f_h(x)=h \lfloor x \rfloor + h (x- \lfloor x \rfloor)=hx,\  g_h(x) = \lfloor hx \rfloor  + (hx - \lfloor hx \rfloor)=hx,    
\end{equation}
so that $f_h(x)=g_h(x)$.

When $h\not= 0$ and $\alpha=0$, we have 
\begin{equation}\label{eqn: alpha=0 inequality}
f_h(x)= \lfloor x \rfloor +  1 ,\  g_h(x) = \lfloor hx \rfloor  + 1,    
\end{equation}
so that $f_h(x) \ge g_h(x)$ when $0<h \le 1$, and $f_h(x) \le g_h(x)$ when $h > 1$.

When $h\not=0$, and $0< \alpha<1$, the set of points where $f_h(x)$ is not differentiable is $x\in \mathbb{Z}^+$, and the set of points where $g_h(x)$ is not differentiable is $hx \in \mathbb{Z}^+$, i.e. $x=k/h$ for $k \in \mathbb{Z}^+$.

Consider the case where $0<h \le 1$ and $0< \alpha <1$, and define
$$F_h(x) := f_h(x)-g_h(x).$$
Since $0< h \le 1$, we have $1/h \ge 1$. 
For any given $k\in \mathbb{Z}^{\ge 0}$, we may assume there exist integers $n_1,n_2,n_3$ (which can potentially be the same integer), such that 
$$\frac{k}{h} \le n_1 \le n_2 \le n_3 < \frac{k+1}{h}.$$
Moreover, without loss of generality, we may also assume $n_1$ is the smallest integer satisfies $k/h \le n_1$, and $n_3$ is the largest integer satisfies $n_3 <(k+1)/h$. 

First, suppose $k/h < x < n_1$, then $$F_h(x)= h^\alpha (n_1 -1) + h^\alpha (x-(n_1-1))^\alpha - k -(hx -k)^\alpha ,$$
so that 
$$F'_h(x)= \alpha h^\alpha(x-n_1+1)^{\alpha-1}-\alpha h (hx-k)^{\alpha-1}.$$
Note that 
\[\begin{array}{crcl}
     & F'_h(x) & \le  & 0 \\
\iff & \alpha h^\alpha(x-n_1+1)^{\alpha-1} &\le & \alpha h (hx-k)^{\alpha-1}  \\
\iff & h^{\alpha-1} (x-n_1+1)^{\alpha-1} &\le & (hx-k)^{\alpha-1}  \\
\iff & (hx-k)^{1-\alpha} &\le & [h(x-n_1+1)]^{1-\alpha}  \\
\iff & hx-k &\le & h(x-n_1+1) \\
\iff & h(n_1-1) &\le & k \\
\iff & n_1-1 &\le &k/h,
\end{array}\]
and the last inequality holds by the assumption above.

Next, consider the case where $n_3<x < (k+1)/h$, then 
$$F_h(x)= h^\alpha n_3 + h^\alpha(x-n_3)^\alpha -k -(hx-k)^\alpha,$$
and 
$$F'_h(x) = \alpha h^\alpha (x-n_3)^{\alpha -1}-\alpha h (hx-k)^{\alpha - 1} .$$
Note that 
\[\begin{array}{crcl}
 & F'_h(x) & \ge  & 0  \\
\iff & \alpha h^\alpha (x-n_3)^{\alpha -1} & \ge & \alpha h (hx-k)^{\alpha - 1} \\
\iff & (hx-k)^{1-\alpha} & \ge & h^{1-\alpha}(x-n_3)^{1-\alpha}  \\
\iff & hx-k & \ge & hx-hn_3 \\
\iff & hn_3 &\ge & k  \\
\iff & n_3 &\ge & k/h,  \\
\end{array}\]
and the last inequality holds by the assumption above.

Lastly, suppose $n_1 < n_3$, then there exists $n_2,n_2 +1$ such that $k/h \le n_2 < n_2+1 < (k+1)/h$, and consider the case where $n_2 < x < n_2 + 1$. 
Then we have 
$$F_h(x) = h^\alpha n_2 + h^\alpha(x-n_2)^\alpha -k -(hx -k)^\alpha ,$$
and 
$$F'_h(x) = \alpha h^\alpha (x-n_2)^{\alpha -1} - \alpha h (hx -k)^{\alpha-1} .$$
Note that 
\[\begin{array}{crcl}
 & F'_h(x) & \ge  & 0  \\
\iff & \alpha h^\alpha (x-n_2)^{\alpha -1} & \ge & \alpha h (hx-k)^{\alpha - 1} \\
\iff & (hx-k)^{1-\alpha} & \ge & h^{1-\alpha}(x-n_2)^{1-\alpha}  \\
\iff & hx-k & \ge & hx-hn_2 \\
\iff & hn_2 &\ge & k  \\
\iff & n_2 &\ge & k/h,  \\
\end{array}\]
and the last inequality holds by the assumption above.

Therefore, for $0 < h \le 1$, $0 < \alpha <1$, $F_h(x)$ has potential local minimum when $x\in \mathbb{Z}^{\ge 0}$, and this is what we shall calculate in the following.

Let $x = N$ for some $N \in \mathbb{Z}^{\ge 0}$, then $$F_h(N) = h^\alpha N  + h^\alpha(N- N)^\alpha - \lfloor hN \rfloor - (hN -\lfloor hN \rfloor)^\alpha.$$
If $N=0$, we have 
$$F_h(N) = h^\alpha\cdot 0 + h^\alpha \cdot 0^\alpha - \lfloor 0 \rfloor -  (0 - \lfloor 0 \rfloor)^\alpha = h^\alpha \cdot 0^\alpha -0^\alpha =0.$$

When $N \ge 1$, since we assume $0< h \le 1$ in the first place, we may further break it up into the union of the following intervals,
$$\frac{k}{N} \le h < \frac{k+1}{N}, \text{where } k\in \mathbb{Z}^{\ge 0},\  k+1 \le N.$$
Let $G(h):=F_h(N)$, then $$G(h) = h^\alpha N -k -(hN -k)^\alpha,$$ 
so that
$$ G'(h)=\alpha N h^{\alpha-1}  - \alpha N (hN -k )^{\alpha-1}, \ 
G''(h) = \alpha(\alpha-1)N h^{\alpha-2}-\alpha(\alpha-1)N^2(hN-k)^{\alpha-2}.$$
If $N=1$, direct calculation gives $F_h(1)=h^\alpha -\lfloor h \rfloor - (h- \lfloor h \rfloor)^\alpha$.
Since  $0< h \le 1$, we have $F_h(1)=0$, and in the following calculation, we may assume $N \ge 2$.

When $G'(h)=0$, we have $h=k/(N-1)$, and 
\begin{eqnarray*}
G''\left(\frac{k}{N-1}\right) 
&=&
\alpha(\alpha-1)N \left[ \left(\frac{k}{N-1} \right)^{\alpha -2}  - N \left( \frac{k N}{N-1} -k\right)^{\alpha-2}  \right]\\
&=& 
\alpha(\alpha-1)N \left[ \left(\frac{k}{N-1} \right)^{\alpha -2}  - N \left( \frac{k}{N-1}\right)^{\alpha-2}  \right] \\
&=& 
\alpha(\alpha-1)N (1-N) \left(\frac{k}{N-1} \right)^{\alpha -2} \\
&\ge& 
0.    
\end{eqnarray*}
This implies $h=k/(N-1)$ is a local minimum for $G(h)$ in the assumed interval $k/N \le h < (k+1)/N$.
In this case, 
$$G\left(\frac{k}{N-1}\right) = \left(\frac{k}{N-1}\right)^{\alpha} N -k - \left(\frac{Nk}{N-1}-k\right)^{\alpha} = (N-1)^{1-\alpha} k^\alpha -k,$$
and $(N-1)^{1-\alpha} k^\alpha -k \ge 0$ is equivalent to $1+k \le N$ which is exactly what we assumed for $k$.
Moreover, $k/(N-1) < (k+1)/N$ if and only if $k+1 <N$, and when $k/(N-1)$ is an endpoint of the assumed  interval, we have the following results.

Finally, when $h=k/N$ or $(k+1)/N$, we have 
$$G\left(\frac{k}{N}\right) = \left(\frac{k}{N}\right)^{\alpha} N -k -\left(\frac{Nk}{N}-k\right)^{\alpha}= N^{1-\alpha}k^\alpha -k\ge 0, $$
since $k < k+1 \le N$, and 
$$G\left(\frac{k+1}{N}\right) = \left(\frac{k+1}{N}\right)^{\alpha} N -k -\left(\frac{N(k+1)}{N}-k\right)^{\alpha}=(k+1)^{\alpha} N^{1-\alpha} - (k+1) \ge 0,$$
since $k+1 \le N$.
Hence, we have $F_h(N) = G(h) \ge 0$, which gives $f_h(x) \ge g_h(x)$ and the inequality result in (\ref{eqn: M_alpha,c scalar multiple}) when $0\le h\le 1$.

Consider the case where $1 < h$ and $0 \le \alpha \le1$.
Note that the special cases where $\alpha =0,1$ are calculated in (\ref{eqn: alpha=0 inequality}) and (\ref{eqn: alpha=1 equality}), and we may assume $0<\alpha <1$ in the following.
Since $1<h$, we have $1/h <1$, and for any given integer $n \in \mathbb{Z}^{\ge 0} $, there exist integers $k_1,k_2,k_3$ (which can potentially be the same integer), such that 
$$n \le \frac{k_1}{h} \le \frac{k_2}{h} \le \frac{k_3}{h} < n+1 .$$
Moreover, we may assume $k_1$ is the smallest integer satisfies $n \le k_1/h$, and $k_3$ is the largest integer satisfies $k_3/h < n+1$.

First, suppose $n < x < k_1/h$, then 
$$F_h(x) = h^\alpha n + h^\alpha(x-n)^\alpha - (k_1-1) - (hx - (k_1-1))^\alpha,$$
so that 
$$F'_h(x) = \alpha h^\alpha (x-n)^{\alpha-1} - \alpha h (hx -(k_1-1))^{\alpha-1}.$$
Note that 
\[\begin{array}{crcl}
     & F'_h(x) & \ge  & 0 \\
\iff & \alpha h^\alpha (x-n)^{\alpha-1} &\ge & \alpha h (hx -(k_1-1))^{\alpha-1} \\
\iff & (hx -(k_1-1))^{1-\alpha} &\ge & h^{1-\alpha} (x-n)^{1-\alpha}  \\
\iff & hx -(k_1-1) &\ge & h(x-n)  \\
\iff & hn &\ge & k_1-1 \\
\iff & n &\ge & (k_1-1)/h, \\
\end{array}\]
and the last inequality holds by the assumption on $k_1$.

Next, suppose $k_3/h <x < n+1$, then 
$$F_h(x) = h^\alpha n + h^\alpha(x-n)^\alpha -k_3-(hx-k_3)^\alpha,$$
so that 
$$F'_h(x)= \alpha h^\alpha (x-n)^{\alpha-1} - \alpha h (hx - k_3)^{\alpha-1}.$$
Note that 
\[\begin{array}{crcl}
     & F'_h(x) & \le  & 0 \\
\iff & \alpha h^\alpha (x-n)^{\alpha-1} &\le & \alpha h (hx - k_3)^{\alpha-1}  \\
\iff & (hx - k_3)^{1-\alpha} &\le & h^{1-\alpha} (x-n)^{1-\alpha}  \\
\iff & hx - k_3 &\le & h(x-n)  \\
\iff & hn &\le & k_3 \\
\iff & n &\le & k_3/h,
\end{array}\]
and the last inequality holds by the assumption above.

Finally, suppose $k_1<k_3$, then there exists $k_2, k_2+1$ such that $n \le k_2/h < (k_2+1)/h < n+1$, and consider the case where $k_2/h <x < (k_2+1)/h $. 
In this case, we have 
$$F_h(x) = h^\alpha n + h^\alpha(x-n)^\alpha -k_2-(hx-k_2)^\alpha,$$
and 
$$F'_h(x)= \alpha h^\alpha (x-n)^{\alpha-1} - \alpha h (hx - k_2)^{\alpha-1}.$$
Note that 
\[\begin{array}{crcl}
     & F'_h(x) & \le  & 0 \\
\iff & \alpha h^\alpha (x-n)^{\alpha-1} &\le & \alpha h (hx - k_2)^{\alpha-1}  \\
\iff & (hx - k_2)^{1-\alpha} &\le & h^{1-\alpha} (x-n)^{1-\alpha}  \\
\iff & hx - k_2 &\le & h(x-n)  \\
\iff & hn &\le & k_2 \\
\iff & n &\le & k_2/h,
\end{array}\]
and the last inequality holds by the assumption above.

Therefore, given $h > 1$ and $0 < \alpha < 1$, $F_h(x)$ has potential local maximum when $x= k/h$ for some $k \in \mathbb{Z}^{\ge 0}$, and we will calculate these values in the following.

Let $N, k \in \mathbb{Z}^{\ge 0}$, and $k \ge N+1$, we may decompose $h>1$ into the (non-disjoint) union of following intervals,
$$\frac{k}{N+1} <h \le \frac{k}{N}, \text{ where } k=N+1,N+2,N+3,\ldots.$$ 

When $N=0$, we have $ k<h \le \infty$, so that $0 \le k/h <1$.  
This gives 
$$
F_h \left( \frac{k}{h} \right) 
= 
h^\alpha \left\lfloor \frac{k}{h} \right\rfloor + h^\alpha \left(  \frac{k}{h} - \left\lfloor \frac{k}{h} \right\rfloor\right)^\alpha - \left\lfloor h\cdot \frac{k}{h} \right\rfloor - \left( h \cdot \frac{k}{h}- \left\lfloor h\cdot \frac{k}{h} \right\rfloor \right)^\alpha    
=
k^\alpha -k \le 0,
$$
since $k \ge 1$ by assumption.

When $N \ge 1$, note that 
$$\frac{k}{N+1} <h \le \frac{k}{N}  \iff N \le \frac{k}{h} < N+1,$$
and this gives 
$$\tilde{G}(h):=F_h \left( \frac{k}{h} \right) = h^\alpha N + h^\alpha \left( \frac{k}{h}-N \right)^\alpha -k -(k-k)^\alpha = h^\alpha N +  (k-hN)^\alpha -k.$$
This gives 
$$\tilde{G}'(h)= \alpha N h^{\alpha -1} - \alpha N (k - hN)^{\alpha-1},$$
and note that 
$$\tilde{G}'(h) < 0 \iff (k-hN)^{1-\alpha} < h^{1-\alpha} \iff k-hN < h \iff \frac{k}{N+1} < h,$$
where the last inequality holds by assumption.
Hence, for $k/(N+1) <h \le k/N$,
\begin{eqnarray*}
\tilde{G}(h) 
&\le& 
\tilde{G}\left( \frac{k}{N+1}\right) \\
&=& \left( \frac{k}{N+1} \right)^\alpha N +\left(k-\frac{kN}{N+1}\right)^\alpha -k \\
&= & 
\left( \frac{k}{N+1} \right)^\alpha (N+1 ) -k \\
&=& 
k^\alpha (N+1)^{1-\alpha} -k \\
&\le& 0,
\end{eqnarray*}
where the last inequality follows from $k \ge N+1$.
Hence, we have $F_h(k/h)=\tilde{G}(h) \le 0$, which gives $f_h(x) \le g_h(x)$. 
When $h=1$, $f_h(x)=g_h(x)$ trivially, so that we have the inequality in (\ref{eqn: M_alpha,c scalar multiple}) when $h \ge 1$.

\end{proof}

When identifying the sum of a transport multi-path as one transport path,  we will see in the following proposition that the $\mathbf{M}_{\alpha,c}$ cost on transport paths potentially has a lower total cost as compared to the cost on transport multi-path.

\begin{proposition}
Let $\vec{T}=(T_1,T_2,\cdots,T_k,\cdots) \in Path_c(\mu^-,\mu^+)$ be a transport multi-path defined as in (\ref{eqn: multi-path}).
For each $k \in \mathbb{N}$, suppose $T_k = \underline{\underline{\tau}}(M_k,\theta_k(x),\xi_k(x))$
such that $\theta_k(x) \le \|\mu_k^-\| = \|\mu_k^+\|$ for all $x \in M_k$,
then 
$$\mathbf{M}_{\alpha,c}(T) \le \mathbf{M}_\alpha(\vec{T}), \text{ for } T = \sum_{k=1}^\infty T_k.$$
Note that for the above inequality, the left side cost is defined using equation (\ref{eqn: cost with capacity}), while the right hand side cost is defined using equation (\ref{eqn: cost with multi-path}).

\end{proposition}

\begin{proof}
For each $k \ge 1$, definition of transport multi-path gives $\|\mu_k^-\|= \|\mu_k^+\| \le c$, which implies $\theta_k(x) \le c$.
Denote $$M_k^0:= \{x \in M_k: \theta_k(x) <c\}, \text{ and }  M_k^1:= \{x \in M_k: \theta_k(x) =c\},$$
then we have 
\begin{eqnarray*}
\mathbf{M}_\alpha(T_k)
&:=& 
\int_{M_k} \theta(x)^\alpha d\mathcal{H}^1  \\
&=& 
\int_{M_k^0} \theta(x)^\alpha d\mathcal{H}^1 + \int_{M_k^1} \theta(x)^\alpha d\mathcal{H}^1    \\
&=&
\int_{M_k^0} c^\alpha \left(\frac{\theta(x)}{c}\right)^\alpha d\mathcal{H}^1 + \int_{M_k^1} c^\alpha \left(\frac{\theta(x)}{c}\right)^\alpha d\mathcal{H}^1  \\
&=&
c^\alpha \cdot \int_{M_k^0} \left\lfloor \frac{\theta(x)}{c}\right\rfloor + \left(\frac{\theta(x)}{c}-\left\lfloor \frac{\theta(x)}{c}\right\rfloor \right)^\alpha d\mathcal{H}^1 
+ 
c^\alpha \cdot \int_{M_k^1}\left\lfloor \frac{\theta(x)}{c}\right\rfloor + \left(\frac{\theta(x)}{c}-\left\lfloor \frac{\theta(x)}{c}\right\rfloor \right)^\alpha d\mathcal{H}^1  \\
&=&
\mathbf{M}_{\alpha,c}(T_k),
\end{eqnarray*}
so that 
$$\mathbf{M}_{\alpha,c}(T) \le \sum_{k=1}^\infty \mathbf{M}_{\alpha,c}(T_k) = \sum_{k=1}^\infty \mathbf{M}_{\alpha}(T_k) = \mathbf{M}_\alpha(\vec{T}). $$

\end{proof}

\begin{remark}
As illustrated in the beginning of this section, given $T = \underline{\underline{\tau}}(M,\theta(x),\xi(x))\in Path(\mu^-,\mu^+)$, the total cost 
$$
\mathbf{M}_{\alpha,c}(T):= c^\alpha\cdot \int_{M} \left\lfloor\frac{\theta(x)}{c}\right\rfloor  + \left( \frac{\theta(x)}{c} - \left\lfloor\frac{\theta(x)}{c}\right\rfloor \right)^\alpha d\mathcal{H}^1,
$$
is dependent on $c$. Here, we will analyze the behavior of the total $\mathbf{M}_{\alpha,c}$ cost when $c$ approaches $0$ or $\infty$.

When $\theta(x) < N$ for $\mathcal{H}^1$ almost every $x \in M$, then $$\lim_{c \to \infty} \left\lfloor\frac{\theta(x)}{c}\right\rfloor =0, \text{ for } \mathcal{H}^1-a.e \ x\in M,$$
and we have 
$$\lim_{c\to \infty} \mathbf{M}_{\alpha,c}(T) = \lim_{c \to \infty} c^\alpha \cdot \int_{M} \left( \frac{\theta(x)}{c}\right)^\alpha d\mathcal{H}^1
= 
\lim_{c \to \infty}  \int_{M} \theta(x)^\alpha d\mathcal{H}^1
=  \int_{M} \theta(x)^\alpha d\mathcal{H}^1.$$
This shows that $\mathbf{M}_{\alpha,c}$ cost equals $\mathbf{M}_\alpha$ cost when $c \to \infty$.

When $\theta(x) >0$ for $\mathcal{H}^1$ almost every $x\in M$, then by definition of $\left\lfloor\theta(x)/c \right\rfloor$, we have 
$$  \frac{\theta(x)}{c} - \left\lfloor\frac{\theta(x)}{c}\right\rfloor  \le 1 .$$
Let $n \in \mathbb{Z}^{+}$, and since $\theta(x) >0$ for $\mathcal{H}^1-a.e. \ x \in M$, we have
$$ \mathcal{H}^1 \left( \left\{x\in M : \theta(x) > \frac{1}{n}\right\}\right) >0,$$
so that 
\begin{eqnarray*}
\lim_{c \to 0} \int_M \left\lfloor\frac{\theta(x)}{c}\right\rfloor  d\mathcal{H}^1 
&\ge& 
\lim_{c \to 0} \int_{\left\{x\in M : \theta(x) > \frac{1}{n}\right\}} \left\lfloor\frac{1/n}{c}\right\rfloor  d\mathcal{H}^1 \\
&\ge&
\int_{\left\{x\in M : \theta(x) > \frac{1}{n}\right\}} \left\lfloor\frac{1/n}{1/n^2}\right\rfloor  d\mathcal{H}^1 \\
&=&
n \cdot \mathcal{H}^1 \left( \left\{x\in M : \theta(x) > \frac{1}{n}\right\}\right),
\end{eqnarray*}
for all $n \in \mathbb{Z}^{>0}$. Therefore, taking $n \to \infty$, we get 
$$\lim_{c \to 0} \int_M \left\lfloor\frac{\theta(x)}{c}\right\rfloor  d\mathcal{H}^1 = \infty.$$

The above calculation gives the following ``estimation'', suppose $\mathcal{H}^1(M) < \infty$, then
$$1 \le \lim_{c \to 0} \frac{ \displaystyle c^\alpha\cdot \int_{M} \left\lfloor\frac{\theta(x)}{c}\right\rfloor  + \left( \frac{\theta(x)}{c} - \left\lfloor\frac{\theta(x)}{c}\right\rfloor \right)^\alpha d\mathcal{H}^1}{ \displaystyle c^\alpha\cdot \int_{M} \left\lfloor\frac{\theta(x)}{c}\right\rfloor  d\mathcal{H}^1} 
\le 1 + \lim_{c \to 0} \frac{\mathcal{H}^1 (M)}{\displaystyle \int_{M} \left\lfloor\frac{\theta(x)}{c}\right\rfloor  d\mathcal{H}^1} =1.
$$
Hence, when $c$ gets sufficiently small as compared to $\theta(x)$ for $x \in M$, the total $\mathbf{M}_{\alpha,c}$ cost  approximately equals to the integer part,
$$c^\alpha \cdot \int_{M} \left\lfloor\frac{\theta(x)}{c}\right\rfloor  d\mathcal{H}^1 .$$

Moreover, given $x \in \mathbb{R}$, we have $ x-1 \le \lfloor x \rfloor \le x$, and this gives
$$c^{\alpha - 1}\theta(x) -c^\alpha   \le 
c^\alpha  \left\lfloor\frac{\theta(x)}{c}\right\rfloor  + c^\alpha  \left( \frac{\theta(x)}{c} - \left\lfloor\frac{\theta(x)}{c}\right\rfloor \right)^\alpha 
\le c^{\alpha-1} \theta(x) + c^\alpha .$$
When $\alpha =1$, 
$$\lim_{c \to 0} c^\alpha  \left\lfloor\frac{\theta(x)}{c}\right\rfloor  + c^\alpha  \left( \frac{\theta(x)}{c} - \left\lfloor\frac{\theta(x)}{c}\right\rfloor \right)^\alpha =\lim_{c \to 0}  c \cdot \frac{\theta(x)}{c} = \theta(x).$$
When $0 \le \alpha <1$, and $\theta(x) >0$,
$$\lim_{c \to 0} c^\alpha  \left\lfloor\frac{\theta(x)}{c}\right\rfloor  + c^\alpha  \left( \frac{\theta(x)}{c} - \left\lfloor\frac{\theta(x)}{c}\right\rfloor \right)^\alpha \ge \lim_{c \to 0} c^{\alpha - 1}\theta(x) -c^\alpha  =\infty.$$

\end{remark}

\begin{remark}\label{rem: 2points to 1 point calculation}
Similar to \cite[Example 2.1]{xia2003} we may also calculate the ``angles'' at an intersection vertex in the simple case when transporting weight from 2 points to 1 point.
Here, the detailed calculation will be conducted by assuming the vertices and transport path are supported in $\mathbb{R}^2$. This is because coordinate calculation will become overwhelmingly complicated when dimension is a general number $n$.

Note that the calculation below only shows aggregation is better than separation, it does not tell us directly what an optimal transport path looks like in the case of 2 points to 1 point.

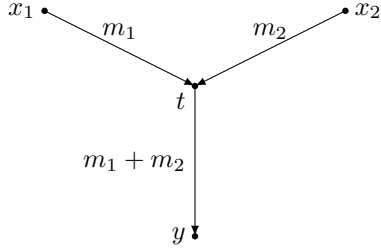
\begin{figure}[h]

\begin{tikzpicture}[>=latex]
\filldraw[black] (-2,3) circle (1pt) node[anchor=east]{$x_1$};
\filldraw[black] (2,3) circle (1pt) node[anchor=west]{$x_2$};
\filldraw[black] (0,0) circle (1pt) node[anchor=east]{$y$};
\filldraw[black] (0,2) circle (1pt) node[anchor = east]{ };
\filldraw[black] (0,1.8) circle (0pt) node[anchor = east]{$t$};

\filldraw[black] (-1, 2.5) circle (0pt) node[anchor = south]{$m_1$};
\filldraw[black] (1, 2.5) circle (0pt) node[anchor = south]{$m_2$};
\filldraw[black] (0, 1) circle (0pt) node[anchor = east]{$m_1 + m_2$};

\draw[->] (-2,3) -- (0,2);
\draw[->] (2,3) -- (0,2);
\draw[->] (0,2) -- (0,0);

\end{tikzpicture}
\caption{Angles at the intersection of transport path $T$.}

\label{fig: angles at intersection }
\end{figure}

Let $x_1 = (x_{1,1},x_{1,2}), \ x_2 = (x_{2,1},x_{2,2}), \ y=(y_1,y_2), \ t= (t_1,t_2)$, and 
$$k_1 = c^\alpha  \left\lfloor\frac{m_1}{c}\right\rfloor  + c^\alpha  \left( \frac{m_1}{c} - \left\lfloor\frac{m_1}{c}\right\rfloor \right)^\alpha, \ 
k_2 = c^\alpha  \left\lfloor\frac{m_2}{c}\right\rfloor  + c^\alpha  \left( \frac{m_2}{c} - \left\lfloor\frac{m_2}{c}\right\rfloor \right)^\alpha, $$  
$$k_3 = c^\alpha  \left\lfloor\frac{m_1+m_2}{c}\right\rfloor  + c^\alpha  \left( \frac{m_1+m_2}{c} - \left\lfloor\frac{m_1+m_2}{c}\right\rfloor \right)^\alpha.$$
The total cost, which we shall express it as a function of $t_1,t_2$, can be formulated as 
$$F(t_1,t_2) = \mathbf{M}_{\alpha,c}(T) = k_1\|x_1 -t\| + k_2 \|x_2 -t\| + k_3 \|y- t\|,$$
where $\|\cdot \|$ is the $2-$norm or Euclidean norm.

Then 
$$\frac{\partial F}{\partial t_1} 
= 
k_1 \frac{t_1 - x_{1,1}}{\|t-x_1\|} + 
k_2 \frac{t_1 - x_{2,1}}{\|t-x_2\|} + 
k_3 \frac{t_1 - y_{1}}{\|t-y\|},$$

$$\frac{\partial F}{\partial t_2} 
= 
k_1 \frac{t_2 - x_{1,2}}{\|t-x_1\|} + 
k_2 \frac{t_2 - x_{2,2}}{\|t-x_2\|} + 
k_3 \frac{t_2 - y_{2}}{\|t-y\|},$$

$$\frac{\partial^2 F}{\partial t_1^2} 
= 
k_1 \frac{(t_2 - x_{1,2})^2}{\|t-x_1\|^3} + 
k_2 \frac{(t_2 - x_{2,2})^2}{\|t-x_2\|^3} + 
k_3 \frac{(t_2 - y_{2})^2}{\|t-y\|^3},$$

$$\frac{\partial^2 F}{\partial t_2^2} 
= 
k_1 \frac{(t_1 - x_{1,1})^2}{\|t-x_1\|^3} + 
k_2 \frac{(t_1 - x_{2,1})^2}{\|t-x_2\|^3} + 
k_3 \frac{(t_1 - y_{1})^2}{\|t-y\|^3},$$

$$\frac{\partial^2 F}{\partial t_1 \partial t_2} 
= 
-k_1 \frac{(t_1 - x_{1,1})(t_2 - x_{1,2})}{\|t-x_1\|^3} - 
k_2 \frac{(t_1 - x_{2,1})(t_2 - x_{2,2})}{\|t-x_2\|^3} - 
k_3 \frac{(t_1 - y_{1})(t_2 - y_{2})}{\|t-y\|^3}.$$

Moreover,
$$\frac{\partial^2 F}{\partial t_1^2} \frac{\partial^2 F}{\partial t_2^2} - \left( \frac{\partial^2 F}{\partial t_1 \partial t_2} \right)^2 $$
consists of ``power'' terms 
$$k_1 \frac{(t_2 - x_{1,2})^2}{\|t-x_1\|^3} \cdot k_1 \frac{(t_1 - x_{1,1})^2}{\|t-x_1\|^3} -     
k_1^2 \frac{(t_1 - x_{1,1})^2(t_2 - x_{1,2})^2}{\|t-x_1\|^6} =0,$$  
and ``cross-product'' terms,
$$k_1 k_2 \frac{(t_2 - x_{1,2})^2(t_1 - x_{2,1})^2+ (t_2 - x_{2,2})^2 (t_1 - x_{1,1})^2}{\|t-x_1\|^3 \|t-x_2\|^3} - 
k_1 k_2  \frac{2(t_1 - x_{1,1})(t_2 - x_{1,2})(t_1 - x_{2,1})(t_2 - x_{2,2})}{\|t-x_1\|^3 \|t-x_2\|^3} \ge 0.$$

When $$\frac{\partial^2 F}{\partial t_1^2} \frac{\partial^2 F}{\partial t_2^2} - \left( \frac{\partial^2 F}{\partial t_1 \partial t_2} \right)^2 =0,$$
using the ``cross-product'' terms, we have the following equations,
$$(t_2 - x_{1,2})(t_1 - x_{2,1}) = (t_2 - x_{2,2}) (t_1 - x_{1,1}) ,$$
$$(t_2 - x_{1,2})(t_1 - y_{1}) = (t_2 - y_{2}) (t_1 - x_{1,1}) ,$$
$$(t_2 - x_{2,2})(t_1 - y_{1}) = (t_2 - y_{2})(t_1 - x_{2,1}) .$$
Without loss of generality, we may translate the transport path such that $y = (y_1,y_2)=(0,0)$. 
Also, we may assume $x_1=(x_{1,1},x_{1,2})$, $x_2=(x_{2,1},x_{2,2})$, and $y=(0,0)$ are not on the same line, which means 
$$x_{1,1}x_{2,2} - x_{1,2}x_{2,1} \not=0 .$$
This is because when all three points are on the same line segment(s), point $t$ equals $x_1$ or $x_2$ or $y$, depending on which weight is merged to which, and the angle  $\angle x_1 t x_2$ equals to either $0$ or $\pi$.

In this case, the three equations above can be simplified into 
$$x_{2,2}t_1 -x_{2,1}t_2 = x_{1,2}t_1 - x_{1,1}t_2 + x_{1,1}x_{2,2}- x_{1,2}x_{2,1} ,  $$
$$x_{1,2}t_1  = x_{1,1} t_2 ,$$
$$x_{2,2}t_1 = x_{2,1}t_2 .$$
Note that $$x_{1,2}x_{2,1}t_2 = x_{2,2} x_{1,2} t_1 =  x_{2,2}x_{1,1} t_2$$ and $x_{1,1}x_{2,2} - x_{1,2}x_{2,1} \not=0$ implies $t_2 =0$.
Similarly, 
$$x_{2,1}x_{1,2}t_1 = x_{2,1}x_{1,1}t_2 = x_{1,1} x_{2,2}t_1$$ 
and $x_{1,1}x_{2,2} - x_{1,2}x_{2,1} \not=0$ implies $t_1 =0$.
Substitute $t_1 =t_2 =0$ into the first equation gives 
$$x_{1,1}x_{2,2}- x_{1,2}x_{2,1} =0 ,$$
contradicting our assumption above.
Therefore, the assumption of $x_1$, $x_2$, $y=(0,0)$ not on the same line segment gives 
$$\frac{\partial^2 F}{\partial t_1^2} \frac{\partial^2 F}{\partial t_2^2} - \left( \frac{\partial^2 F}{\partial t_1 \partial t_2} \right)^2 >0 $$
Moreover, this assumption also gives 
$$\frac{\partial^2 F}{\partial t_1^2} > 0.$$
Finally, the existence of $t'=(t_1',t_2')$ such that 
$$\frac{\partial F}{\partial t_1} (t'_1,t'_2) =0 = \frac{\partial F}{\partial t_2} (t'_1,t'_2)  ,$$
is equivalent to the existence of a triangle 
\begin{equation}\label{eqn: triangle vector for derivative}
k_1 \vec{n}_1 + k_2 \vec{n}_2+ k_3\vec{n}_3 =0,   
\end{equation}
where $\vec{n}_1, \vec{n}_2,\vec{n}_3,$ are the unit vectors,
$$\vec{n}_1 = \frac{t'-x_1}{\|t' -x_1\|}, \ 
\vec{n}_2 = \frac{t'-x_2}{\|t' -x_2\|},\ 
\vec{n}_3 = \frac{t'-y}{\|t' -y\|}. $$
By using similar calculation as in \ref{eqn: M_alpha,c subadditivity}, we have $k_3 \le k_1 + k_2$, which gives the existence of triangle.
Hence, $t'$ is a local minimum for $F(t) = \mathbf{M}_{\alpha,c}(T)$.
This shows that ``Y'' shaped transport paths has lower total transport cost than ``V'' shaped transport path.

Again, using formula \ref{eqn: triangle vector for derivative} and law of cosine, we can also calculate the angle $\theta_1$ between $-\vec{n}_1$ and $\vec{n}_3$, angle $\theta_2$ between $-\vec{n}_2$ and $\vec{n}_3$, angle $\theta_3$ between $-\vec{n}_1$ and $-\vec{n}_2$.
Direct calculation gives 
$$cos(\theta_1) = \frac{k_1^2 +k_3^2 -k_2^2}{2k_1 k_3} = \frac{(k_1/k_3)^2 + 1 - (k_2/k_3)^2}{2 k_1/k_3},$$
$$cos(\theta_2) = \frac{k_2^2 + k_3^2 -k_1^2}{2k_2k_3} = \frac{(k_2/k_3)^2 + 1 - (k_1/k_3)^2}{2 k_2/k_3}.$$
By parallelogram law we have 
$$\|-k_3 \vec{n}_3\|^2 + \|k_1 \vec{n}_1 - k_2\vec{n}_2\|^2 = 2 \|k_1 \vec{n}_1\|^2 + 2 \|k_2 \vec{n}_2\|^2 ,$$
so that $\|k_1 \vec{n}_1 - k_2\vec{n}_2\|^2 = 2k_1^2 +2k_2^2 - k_3^2$,
which implies 
$$cos(\theta_3) = \cos(\theta_1 + \theta_2) = \frac{k_1^2 +k_2^2 -(2k_1^2 +2k_2^2 - k_3^2)}{2k_1k_2} = \frac{k_3^2 -k_1^2 -k_2^2}{2k_1k_2} = \frac{1-(k_1/k_3)^2-(k_2/k_3)^2}{2(k_1/k_3)(k_2/k_3)}.$$
By the above calculation, the angles are dependent on the ratio of the cost before and after aggregation.

When $c$ gets sufficiently large, i.e. $c > m_1 + m_2$,
we have
$$\left\lfloor \frac{m_1}{c} \right\rfloor = \left\lfloor \frac{m_2}{c} \right\rfloor =\left\lfloor \frac{m_1+ m_2}{c} \right\rfloor =0 .$$
This implies 
$$\lim_{c \to \infty} \frac{k_1}{k_3} =  
\lim_{c \to \infty} \frac{\left\lfloor\frac{m_1}{c}\right\rfloor  + \left( \frac{m_1}{c} - \left\lfloor\frac{m_1}{c}\right\rfloor \right)^\alpha}{\left\lfloor\frac{m_1+m_2}{c}\right\rfloor + \left( \frac{m_1+m_2}{c} - \left\lfloor\frac{m_1+m_2}{c}\right\rfloor \right)^\alpha} = 
\frac{m_1^\alpha}{(m_1+m_2)^\alpha},
$$
and 
$$\lim_{c \to \infty} \frac{k_2}{k_3} =  
\lim_{c \to \infty} \frac{\left\lfloor\frac{m_2}{c}\right\rfloor  + \left( \frac{m_2}{c} - \left\lfloor\frac{m_2}{c}\right\rfloor \right)^\alpha}{\left\lfloor\frac{m_1+m_2}{c}\right\rfloor + \left( \frac{m_1+m_2}{c} - \left\lfloor\frac{m_1+m_2}{c}\right\rfloor \right)^\alpha} = 
\frac{m_2^\alpha}{(m_1+m_2)^\alpha}. $$
Therefore, as $c \to \infty$, $\cos(\theta_1), \cos(\theta_2)$ behaves the same as in the $\mathbf{M}_{\alpha}$ cost assumption

Next, given any positive integer $N \in \mathbb{Z}^{+}$, when $c$ gets sufficiently small, we may assume 
$$\left\lfloor \frac{m_1}{c} \right\rfloor, \left\lfloor \frac{m_2}{c} \right\rfloor, \left\lfloor \frac{m_1+ m_2}{c} \right\rfloor > N  .$$
For each $x \in \mathbb{R}$, we have $ x-1 \le \lfloor x \rfloor \le x$, and this gives the following inequality.
$$
\frac{\frac{m_1}{c} -1}{\frac{m_1+m_2}{c} + 1} \le
\frac{\left\lfloor\frac{m_1}{c}\right\rfloor}{\left\lfloor\frac{m_1+m_2}{c}\right\rfloor +1} \le 
\frac{k_1}{k_3}
\le \frac{\left\lfloor\frac{m_1}{c}\right\rfloor + 1}{\left\lfloor\frac{m_1+m_2}{c}\right\rfloor}
\le \frac{\frac{m_1}{c} + 1}{\frac{m_1+m_2}{c} -1}.$$
Since
$$\lim_{c\to 0} \frac{\frac{m_1}{c} -1}{\frac{m_1+m_2}{c} + 1} = \frac{m_1}{m_1+m_2} = \lim_{c \to 0}\frac{\frac{m_1}{c} + 1}{\frac{m_1+m_2}{c} -1},$$
we have $$\lim_{c \to 0} \frac{k_1}{k_3} = \frac{m_1}{m_1+m_2},
\text{ and similarly } 
\lim_{c \to 0} \frac{k_2}{k_3} = \frac{m_2}{m_1+m_2}.$$
In this case, we have 
$$\lim_{c \to 0} \cos(\theta_1) = \lim_{c \to 0} \cos(\theta_2) = \lim_{c \to 0} \cos(\theta_1 + \theta_2) =1,$$
which implies $x_1,x_2,t',y$ are on the same line segment.
\end{remark}

\section{Cycle-free property of transport paths.}
Motivated by \cite[Proposition 2.1]{xia2003}, which gives optimal transport paths under $\mathbf{M}_\alpha$ cost are cycle-free. 
We would like to investigate similar results when using $\mathbf{M}_{\alpha,c}$ cost.
Note that the notion of \textit{acyclic} defined using subcurrent given in (\ref{eqn: subcurrent}) is orientation ``sensitive'', while the notion of \textit{cycle-free} given in Definition \ref{def: cycle-free current} is orientation ``insensitive''.

\begin{proposition}\label{prop: integer and decimal acyclic components}

Given $T = \underline{\underline{\tau}}(M,\theta(x),\xi(x)) \in Path (\mu^-, \mu^+), \alpha \in [0,1]$, $c>0$. 
Suppose $T$ is cyclic and
$$\min_{x\in M} \left\lfloor \frac{\theta(x)}{c} \right\rfloor  \ge 1,$$
then we can find a transport path $T_0 =\underline{\underline{\tau}}(M_0,\theta_0(x),\xi(x)) \in Path(\mu^-,\mu^+)$, such that 
$$\mathbf{M}_{\alpha,c}(T_0) \le  \mathbf{M}_{\alpha,c}(T), M_0 \subseteq M,\ \text{and}  \min_{x\in M_0} \left\lfloor \frac{\theta_0(x)}{c} \right\rfloor  =0 .$$
\end{proposition}

\begin{proof}
Suppose $T$ is cyclic, then there exists $S=\underline{\underline{\tau}}(N,\phi(x),\zeta(x))$ on $T$, such that $\partial S =0$. Hence, we have 
$$\mathcal{H}^1(N\backslash M) =0, \text{ and }  \zeta(x) = \pm \xi(x), \text{ for all } x \in N. $$
Denote $R:= \underline{\underline{\tau}}(N, 1 ,\zeta(x))$, $t \cdot R:=\underline{\underline{\tau}}(N, t ,\zeta(x))$ for $t \in \mathbb{R}$, and
$$N_1:= \{x\in N: \zeta(x)=\xi(x)\}, \ N_2:= \{x\in N: \zeta(x)=-\xi(x)\}.$$
Here, without loss of generality, we may assume $\mathcal{H}^1(N_1)-\mathcal{H}^1(N_2) \le 0$. Otherwise, we may simply switching the indexes to let the inequality hold.

Let $$n_0 = \min_{x\in N_2} \left\lfloor \frac{\theta(x)}{c} \right\rfloor \ge \min_{x\in M} \left\lfloor \frac{\theta(x)}{c} \right\rfloor \ge 1,$$
and consider
$$T_0:=T+n_0c \cdot R = \underline{\underline{\tau}}(N_1,\theta(x) + n_0c,\xi(x) )+ \underline{\underline{\tau}}(N_2,\theta(x) - n_0c,\xi(x) )+ \underline{\underline{\tau}}(M\backslash (N_1 \cup N_2),\theta(x),\xi(x) ).$$
Suppose $n \in \mathbb{Z}$ and $0 \le n \le n_0$, then  
\begin{eqnarray*}
F(nc) 
&:=& 
\mathbf{M}_{\alpha,c}(T+ nc \cdot R)-\mathbf{M}_{\alpha,c}(T)  \\ 
&=& \hphantom{-}
c^\alpha \int_{N_1} \left\lfloor\frac{\theta(x)+nc}{c}\right\rfloor  + \left( \frac{\theta(x)+nc}{c} - \left\lfloor\frac{\theta(x)+nc}{c}\right\rfloor \right)^\alpha d\mathcal{H}^1 \\
& & -  
c^\alpha \int_{N_1} \left\lfloor\frac{\theta(x)}{c}\right\rfloor  + \left( \frac{\theta(x)}{c} 
- \left\lfloor\frac{\theta(x)}{c}\right\rfloor \right)^\alpha d\mathcal{H}^1  \\
& & +
c^\alpha \int_{N_2} \left\lfloor\frac{\theta(x)-nc}{c}\right\rfloor  + \left( \frac{\theta(x)-nc}{c} - \left\lfloor\frac{\theta(x)-nc}{c}\right\rfloor \right)^\alpha d\mathcal{H}^1 \\
& & -  
c^\alpha \int_{N_2} \left\lfloor\frac{\theta(x)}{c}\right\rfloor  + \left( \frac{\theta(x)}{c} 
- \left\lfloor\frac{\theta(x)}{c}\right\rfloor \right)^\alpha d\mathcal{H}^1  \\
&=& \hphantom{+}
c^\alpha \int_{N_1} \left\lfloor\frac{\theta(x)}{c}\right\rfloor +n  + \left( \frac{\theta(x)}{c} - \left\lfloor\frac{\theta(x)}{c}\right\rfloor \right)^\alpha d\mathcal{H}^1 -  
c^\alpha \int_{N_1} \left\lfloor\frac{\theta(x)}{c}\right\rfloor  + \left( \frac{\theta(x)}{c} 
- \left\lfloor\frac{\theta(x)}{c}\right\rfloor \right)^\alpha d\mathcal{H}^1  \\
& & +
c^\alpha \int_{N_2} \left\lfloor\frac{\theta(x)}{c}\right\rfloor -n + \left( \frac{\theta(x)}{c} - \left\lfloor\frac{\theta(x)}{c}\right\rfloor \right)^\alpha d\mathcal{H}^1 
-  
c^\alpha \int_{N_2} \left\lfloor\frac{\theta(x)}{c}\right\rfloor  + \left( \frac{\theta(x)}{c} 
- \left\lfloor\frac{\theta(x)}{c}\right\rfloor \right)^\alpha d\mathcal{H}^1  \\
&=& \hphantom{+}
c^\alpha \int_{N_1} n d\mathcal{H}^1-c^\alpha \int_{N_2} n d\mathcal{H}^1 \\
&=& \hphantom{+}
c^\alpha \cdot n \cdot \left(\mathcal{H}^1(N_1)-\mathcal{H}^1(N_2)\right) \le 0.
\end{eqnarray*}
The above calculation shows that 
$$\mathbf{M}_{\alpha,c}(T_0) =  \mathbf{M}_{\alpha,c}(T+n_0c \cdot R)  \le  \mathbf{M}_{\alpha,c}(T).$$

Also, for $x \in N_2 \subseteq M$, $\theta(x) -n_0 c \ge 0 $, so that $T+n_0c \cdot R$ is still an admissible transport path, with 
$$\partial (T+n_0c \cdot R)= \partial T + n_0c \cdot \partial R = \partial T,$$
so that $T+n_0c \cdot R \in Path(\mu^-,\mu^+)$.
Moreover, suppose $T_0 =\underline{\underline{\tau}}(M_0,\theta_0(x),\xi(x))$ then $M_0 \subseteq M$ and
$$\min_{x\in M_0} \left\lfloor  \frac{\theta_0(x)}{c} \right\rfloor 
= \min_{x\in N_2} \left\lfloor  \frac{\theta_0(x)}{c} \right\rfloor =\min_{x\in N_2} \left\lfloor  \frac{\theta(x)-n_0c}{c} \right\rfloor 
=\min_{x\in N_2} \left\lfloor  \frac{\theta(x)}{c} \right\rfloor -n_0
= 0.$$

\end{proof}

\begin{corollary}\label{cor: interger multiple multiplicity}

Given $T = \underline{\underline{\tau}}(M,\theta(x),\xi(x)) \in Path (\mu^-, \mu^+), \alpha \in [0,1]$, $c>0$. 
Suppose 
$$\min_{x\in M} \theta(x) = n_0 \cdot c, \text{ where } n_0 \in \mathbb{Z}^{\ge 0}. $$
Then we can find a transport path $T_0 =\underline{\underline{\tau}}(M_0,\theta_0(x),\xi(x)) \in Path(\mu^-,\mu^+)$, such that 
$$\mathbf{M}_{\alpha,c}(T_0) \le  \mathbf{M}_{\alpha,c}(T), M_0 \subseteq M,\ \text{and}  \min_{x\in M_0} \theta_0(x) =0 .$$
That being said, $T_0$ is cycle-free.
\end{corollary}

\begin{proof}
If $n_0 =0$, then let $T_0:=T$ gives the desired result.

If $n_0 \ge 1$, consider the transport path $T_0:= T+ n_0c\cdot R$ defined above, which gives that 
$$0 \le \min_{x\in M_0} \theta_0(x) \le \min_{x\in N_2} \theta(x) -n_0\cdot c =0,$$
and we get the desired result.
\end{proof}

\begin{theorem}\label{thm: integer part and decimal part}
Given $T = \underline{\underline{\tau}}(M,\theta(x),\xi(x)) \in Path (\mu^-, \mu^+), \alpha \in [0,1]$, $c>0$. 
Let $R$ be an arbitrary cycle on $T$, such that $R = \underline{\underline{\tau}}(N,\phi(x),\zeta(x))$.
Suppose for any point $x \in N$, we have 
\begin{equation}\label{eqn: decimal sum less than 1}
\max_{x\in N} \left( \frac{\theta(x)}{c} -  \left\lfloor \frac{\theta(x)}{c} \right \rfloor \right)  + \min_{x\in N} \left( \frac{\theta(x)}{c} -  \left\lfloor \frac{\theta(x)}{c} \right \rfloor \right)  \le 1,    
\end{equation}
then we may find $$T_1=\underline{\underline{\tau}}(M,\theta_1(x),\xi(x)), \text{ and } T_2=\underline{\underline{\tau}}(M,\theta_2(x),\xi(x))$$ such that  
$\partial T= \partial ( T_1 + T_2 )$,
$$\theta_1(x) = c \cdot n(x), \text{ where } n(x) \in \mathbb{Z}^+,\  \theta_2(x) \le c.$$
Moreover, both 
$T_1$, $T_2$ are cycle-free transport paths, such that 
$$\mathbf{M}_{\alpha,c}(T_1+T_2) = \mathbf{M}_{\alpha,c}(T_1)+\mathbf{M}_{\alpha,c}(T_2)\le \mathbf{M}_{\alpha,c}(T).$$
\end{theorem}

\begin{proof}
Let 
$$T'_1:= c \cdot \underline{\underline{\tau}} \left(M,\left\lfloor \frac{\theta(x)}{c} \right\rfloor,\xi(x)\right) = \underline{\underline{\tau}} \left(M,c \left\lfloor \frac{\theta(x)}{c} \right\rfloor,\xi(x)\right) ,$$
and denote 
$$T_2':= T-T'_1 = \underline{\underline{\tau}}\left(M,\theta(x)- c \left\lfloor \frac{\theta(x)}{c} \right\rfloor,\xi(x) \right).$$

Suppose $T'_1$ is cycle-free, then setting $T_1=T'_1$.
Suppose there exists a cycle $S_1=\underline{\underline{\tau}}(N_1,\phi_1(x),\zeta_1(x))$ on $T'_1$, 
and let $R_1 := \underline{\underline{\tau}}(N_1, 1 ,\zeta_1(x))$.
Let
$$M'_1 =\left\{ x\in M : \left\lfloor \frac{\theta(x)}{c} \right\rfloor \ge 1 \right\}, \text{ and } n_1 = \min_{x\in M_1'} \left\lfloor \frac{\theta(x)}{c} \right\rfloor \ge 1 .$$
If $M_1' = \emptyset $, then $\left\lfloor \theta(x)/c  \right\rfloor =0$ for all $x \in M$, and $T_1' = \underline{\underline{\tau}}\left(M, 0,\xi(x) \right)$. 
Since $S_1$ is on $T_1'$, we have $\phi_1(x) = 0$, which implies $T_1'$ is cycle-free. 
Therefore, we may assume $M_1' \not= \emptyset$.

By Proposition \ref{prop: integer and decimal acyclic components} or Corollary \ref{cor: interger multiple multiplicity}, we can find 
$T_1 := T'_1 + n_1 c \cdot R_1 = \underline{\underline{\tau}}(M,\theta_1(x),\xi(x))$, such that $$\partial T_1 = \partial T'_1,\  \mathbf{M}_{\alpha,c}(T_1) \le \mathbf{M}_{\alpha,c}(T'_1), \ M_1 \subseteq M'_1,$$
and $T_1$ is cycle-free.
Also, note that 
$$\theta_1(x)= c \cdot \left( \left\lfloor \frac{\theta(x)}{c} \right\rfloor + n_1 \right), \text{ and } \left\lfloor \frac{\theta(x)}{c} \right\rfloor + n_1 \in \mathbb{Z}^+.$$

Similarly, if $T'_2$ is cycle-free, then setting $T_2=T'_2$. 
Suppose there is a cycle $S_2=\underline{\underline{\tau}}(N_2,\phi_2(x),\zeta_2(x))$ on $T'_2$, and let $R_2 := \underline{\underline{\tau}}(N_2, 1 ,\zeta_2(x))$, 
$$n_2 = \min_{x\in N_2}\left( \theta(x)- c \left\lfloor \frac{\theta(x)}{c} \right\rfloor 
\right).$$
Since $S_2$ is on $T_2'$, we have $0 < \phi_2(x) \le \theta(x)- c \lfloor \theta(x)/c \rfloor$ when $x \in N_2$, which further implies $n_2 > 0$.

Similar to the argument as in  \cite[Proposition 2.1]{xia2003}, we may first denote 
$$N_2^+=\{x\in N_2 : \zeta_2(x)= \xi(x) \},\ N_2^-=\{x\in N_2 : \zeta_2(x)= -\xi(x) \}.$$
Let $t \in [0,1]$, then 
\begin{eqnarray*}
T'_2+ tn_2\cdot R_2 
&=&
\underline{\underline{\tau}}\left(N_2^+, \theta(x)- c \left\lfloor \frac{\theta(x)}{c} \right\rfloor + tn_2 ,\xi(x)\right)    + \\
&&
\underline{\underline{\tau}}\left(N_2^-, \theta(x)- c \left\lfloor \frac{\theta(x)}{c} \right\rfloor - tn_2 ,\xi(x)\right)    + \\
&&
\underline{\underline{\tau}}\left(M\backslash N_2, \theta(x)- c \left\lfloor \frac{\theta(x)}{c} \right\rfloor ,\xi(x)\right).
\end{eqnarray*}
Since $R_2$ is a cycle on $T_2$, so that it is also a cycle on $T$, by the assumption in (\ref{eqn: decimal sum less than 1}) we have
\begin{eqnarray*}
\theta(x)- c \left\lfloor \frac{\theta(x)}{c} \right\rfloor \pm tn_2 
&\le& 
\max_{x\in N_2} \left(\theta(x)- c \left\lfloor \frac{\theta(x)}{c} \right\rfloor \right) + t \min_{x\in N_2}\left( \theta(x)- c \left\lfloor \frac{\theta(x)}{c} \right\rfloor \right) \\
&\le &
c\cdot \max_{x\in N_2} \left(\frac{\theta(x)}{c}- \left\lfloor \frac{\theta(x)}{c} \right\rfloor \right) + c\cdot  \min_{x\in N_2}\left( \frac{\theta(x)}{c} - \left\lfloor \frac{\theta(x)}{c} \right\rfloor \right) \\
&\le&
c.
\end{eqnarray*}
Therefore, direct calculation gives
\begin{eqnarray*}
F(t) 
&:=&
\mathbf{M}_{\alpha,c}(T'_2+ tn_2\cdot R_2)-\mathbf{M}_{\alpha,c}(T'_2) \\
&=& 
\int_{N_2^+} \left( \theta(x)- c \left\lfloor \frac{\theta(x)}{c} \right\rfloor +tn_2 \right)^\alpha - \left( \theta(x)- c \left\lfloor \frac{\theta(x)}{c} \right\rfloor \right)^\alpha d\mathcal{H}^1 + \\
&& 
\int_{N_2^-} \left( \theta(x)- c \left\lfloor \frac{\theta(x)}{c} \right\rfloor -tn_2 \right)^\alpha - \left( \theta(x)- c \left\lfloor \frac{\theta(x)}{c} \right\rfloor \right)^\alpha d\mathcal{H}^1,
\end{eqnarray*}

$$F'(t)=\alpha n_2 \int_{N_2^+} \left( \theta(x)- c \left\lfloor \frac{\theta(x)}{c} \right\rfloor +tn_2 \right)^{\alpha-1} d\mathcal{H}^1 
- 
\alpha n_2 \int_{N_2^-} \left( \theta(x)- c \left\lfloor \frac{\theta(x)}{c} \right\rfloor -tn_2 \right)^{\alpha-1} d\mathcal{H}^1,$$
$$F''(t)=\alpha(\alpha-1) n^2_2 \int_{N_2^+} \left( \theta(x)- c \left\lfloor \frac{\theta(x)}{c} \right\rfloor +tn_2 \right)^{\alpha-2} \!d\mathcal{H}^1 
+
\alpha(\alpha-1) n^2_2 \int_{N_2^-} \left( \theta(x)- c \left\lfloor \frac{\theta(x)}{c} \right\rfloor -tn_2 \right)^{\alpha-2} \!d\mathcal{H}^1.$$
By picking an reverse orientation of $R_2$ if necessary, i.e. replace $\zeta_2(x)$ by $-\zeta_2(x)$, we may assume $F'(0) \le 0$.
$F''(t) \le 0$ implies $F'(t)$ is decreasing, and $F'(0) \le 0$ gives minimum occurs when $t=1$, so that $$\mathbf{M}_{\alpha,c}(T'_2+ n_2\cdot R_2)-\mathbf{M}_{\alpha,c}(T'_2) = F(1) \le F(0) =0 .$$

Hence, by defining $T_2:= T'_2 + n_2 \cdot R_2 = \underline{\underline{\tau}}(M,\theta_2(x),\xi(x))$ we have $\partial T_2 = \partial T'_2$,
$$\theta_2(x) = \theta(x)- c \left\lfloor \frac{\theta(x)}{c} \right\rfloor \pm tn_2 \le c \text{ for } t \in [0,1] .$$
Since $$\min_{x \in N_2} \theta_2(x) = \min_{x \in N_2}\left( \theta(x)- c \left\lfloor \frac{\theta(x)}{c} \right\rfloor \right)  - n_2 = \min_{x \in N_2}\left( \theta(x)- c \left\lfloor \frac{\theta(x)}{c} \right\rfloor \right) -\min_{x \in N_2}\left( \theta(x)- c \left\lfloor \frac{\theta(x)}{c} \right\rfloor \right) =0,$$
so that $R_2$ is no longer a cycle on $T_2$. Since we may do this ``cycle reduction'' for arbitrary cycles on $T'_2$, we may ultimately assume that $T_2$ is cycle-free. Also $\partial (T_1 + T_2) = \partial T_1 + \partial T_2 = \partial T'_1 + \partial T'_2= \partial(T'_1+T'_2) = \partial T$.

Lastly, we shall prove 
$$\mathbf{M}_{\alpha,c}(T_1+T_2) = \mathbf{M}_{\alpha,c}(T_1)+\mathbf{M}_{\alpha,c}(T_2)\le \mathbf{M}_{\alpha,c}(T).$$

Note that $$T_1=  \underline{\underline{\tau}}(M,\theta_1(x),\xi(x)), \ \theta_1(x) =  c \cdot \left( \left\lfloor \frac{\theta(x)}{c} \right\rfloor + n_1 \right), \ \left\lfloor \frac{\theta(x)}{c} \right\rfloor + n_1 \in \mathbb{Z}^+,$$
and when $T'_1$ is cycle-free, by Proposition \ref{prop: integer and decimal acyclic components} we may assume $n_1=0$.
This gives 
\begin{eqnarray*}
\mathbf{M}_{\alpha,c}(T_1)+ \mathbf{M}_{\alpha,c}(T_2) 
&=&
c^\alpha \int_{M} \left\lfloor \frac{\theta_1(x)}{c} \right\rfloor   d\mathcal{H}^1 
+
c^\alpha \int_{M} \left( \frac{\theta_2(x)}{c} \right)^\alpha  d\mathcal{H}^1 \\
&=&
c^\alpha \int_M \left\lfloor \frac{\theta_1(x) + \theta_2(x)}{c} \right\rfloor  + \left( \frac{\theta_1(x) + \theta_2(x)}{c} -  \left\lfloor \frac{\theta_1(x) + \theta_2(x)}{c} \right\rfloor \right)^\alpha  d\mathcal{H}^1 \\
&=&
\mathbf{M}_{\alpha,c}(T_1 + T_2).
\end{eqnarray*}

Similarly, we also have
\begin{eqnarray*}
\mathbf{M}_{\alpha,c}(T'_1)+ \mathbf{M}_{\alpha,c}(T'_2) 
&=&
c^\alpha \int_{M} \left\lfloor \frac{\theta(x)}{c} \right\rfloor d\mathcal{H}^1 +
c^\alpha \int_{M} \left( \frac{\theta(x)}{c} - \left\lfloor \frac{\theta(x)}{c} \right\rfloor \right)^\alpha d\mathcal{H}^1 \\
&=&
\mathbf{M}_{\alpha,c}(T'_1+T'_2).
\end{eqnarray*}
Hence, we have
$$\mathbf{M}_{\alpha,c}(T_1+T_2) = \mathbf{M}_{\alpha,c}(T_1)+\mathbf{M}_{\alpha,c}(T_2)\le \mathbf{M}_{\alpha,c}(T'_1)+ \mathbf{M}_{\alpha,c}(T'_2) = \mathbf{M}_{\alpha,c}(T'_1+T'_2) = \mathbf{M}_{\alpha,c}(T). $$
\end{proof}

When $T = \underline{\underline{\tau}}(M,\theta(x),\xi(x))$ and  $T_0 =\underline{\underline{\tau}}(M_0,\theta_0(x),\xi(x))$, the conditions $$ M_0 \subseteq M  \text{ and } \min_{x\in M_0} \left\lfloor \frac{\theta_0(x)}{c} \right\rfloor  =0 ,$$
do not automatically imply $T_0$ is cycle-free or cycle-free is a potential property of an optimal transport path.
This can be demonstrated by the following example.

\begin{example}
Suppose we have a transport path (locally) of the form $$T= m_1 \llbracket \overline{x_1 x_2} \rrbracket + m_2 \llbracket \overline{x_2 x_3} \rrbracket + m_3 \llbracket \overline{x_1 x_3} \rrbracket,$$ where $m \llbracket \overline{x_i x_j} \rrbracket$ represents the rectifiable currents supported on the line segment $\overline{x_ix_j}$, with density $m$, in the direction of $\overrightarrow{x_ix_j}$.
Let $R$ be the counterclockwise defined transport path with constant density $1$, i.e. $R = \llbracket \overline{x_1 x_3} \rrbracket + \llbracket \overline{x_3 x_2} \rrbracket + \llbracket \overline{x_2 x_1} \rrbracket$.

\begin{figure}[h]
\centering
\begin{tikzpicture}[>=latex]

\filldraw[black] (0,0) circle (1pt) node[anchor=east]{$x_1$};
\filldraw[black] (0.25,1.837) circle (1pt) node[anchor=east]{$x_2$};
\filldraw[black] (2.25,1.837) circle (1pt) node[anchor=west]{$x_3$};

\filldraw[black] (0.12, 0.9) circle (0pt) node[anchor = east]{$m_1$};
\filldraw[black] (1.125, 1.8) circle (0pt) node[anchor = south]{$m_2$};
\filldraw[black] (1.125, 0.9) circle (0pt) node[anchor = west]{$m_3$};

\draw[->] (0,0) -- (0.25,1.837);
\draw[->] (0.25,1.837) -- (2.25,1.837);
\draw[->] (0,0) -- (2.25,1.837);

\end{tikzpicture}
\caption{Cyclic transport path $T$.}
\label{fig: cyclic transport path}

\end{figure}
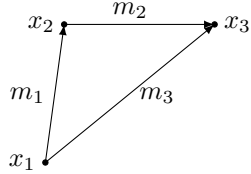

Also, we further assume $\alpha =0.5,\  c=1$, $$m_1 =0.5,\  m_2=1,\ m_3=1, \text{ and } |x_2 -x_1| = 1,\  |x_3 -x_2| = 1,\  |x_3 -x_1| = 1.5.$$
In this specific example, for arbitrary $t \in \mathbb{R}$,  we identify a rectifiable current as follows,
$$t \cdot \underline{\underline{\tau}}(M,\theta(x),\xi(x)) 
= \underline{\underline{\tau}}(M,t\cdot \theta(x),\xi(x))
= \underline{\underline{\tau}}(M,(-t)\cdot \theta(x),-\xi(x)).$$ 
So that 
$t\llbracket \overline{x_i x_j} \rrbracket = (-t) \llbracket\overline{x_j x_i}\rrbracket$
and 
$$tR= (-t)\llbracket \overline{x_3 x_1 } \rrbracket + (-t)\llbracket \overline{x_2 x_3 } \rrbracket + (-t)\llbracket \overline{x_1x_2 } \rrbracket ,$$
which is still a transport paths with potentially negative density values and reversed orientation.

In this case, for all $t\in \mathbb{R}$, we have 
$$T+ tR = (m_1 - t) \llbracket \overline{x_1 x_2} \rrbracket + (m_2 - t) \llbracket \overline{x_2 x_3} \rrbracket + (m_3 + t) \llbracket \overline{x_1 x_3} \rrbracket $$
and 
\begin{eqnarray*}
F(t):=\mathbf{M}_{\alpha,c}(T+ tR) 
&=&
\left\lfloor |m_1-t| \right\rfloor + (|m_1 - t| -\left\lfloor |m_1-t| \right\rfloor)^\alpha |x_2 -x_1| +  \\
&&
\left\lfloor |m_2-t| \right\rfloor + (|m_2 - t| -\left\lfloor |m_2-t| \right\rfloor)^\alpha |x_3 -x_2| +  \\
&&
\left\lfloor |m_3+t| \right\rfloor + (|m_3 +t| -\left\lfloor |m_3+t| \right\rfloor)^\alpha |x_3 -x_1|  .
\end{eqnarray*}
For arbitrary $N \in \mathbb{Z}$, consider $t \in (N,N+1/2) \cup (N + 1/2, N+1)$.

When $N \ge 1$, and $t \in (N,N+1/2)$, 
$$F(t)=[N-1 + (t- \frac{1}{2} - (N-1) )^\alpha] + [N-1 + (t-1 -(N-1))^\alpha] + [N+1 + (t+1 - (N+1))^\alpha]\frac{3}{2}.$$
Then, 
$$\lim_{t \to (N+ 1/2)^-} F(t) = [(N-1) + 1 ^\alpha] + [N-1  + \left( \frac{1}{2}\right)^\alpha] + [N+1 + \left( \frac{1}{2}\right)^\alpha] \frac{3}{2}, $$
and 
$$F'(t) = \alpha (t- N + \frac{1}{2})^{\alpha - 1} + \alpha (t-N)^{\alpha -1} + \frac{3}{2}\alpha (t-N)^{\alpha -1} > 0 .$$

When $N \ge 1$, and $t \in (N+ 1/2, N+1)$,
$$F(t)= [N + (t- \frac{1}{2} - N)^\alpha] + [N-1 + (t-1 - (N-1))^\alpha] + [N+1 + (t+1- (N+1)^\alpha)]\frac{3}{2} .$$
Then,
$$\lim_{t \to (N + 1/2)^+} F(t) = [N + 0^\alpha] + [N-1 + \left( \frac{1}{2} \right)^\alpha]+ [N+1 + \left( \frac{1}{2}\right)^\alpha]\frac{3}{2},$$
and 
$$F'(t) = \alpha (t -\frac{1}{2}-N)^{\alpha -1}  + \alpha (t-N)^{\alpha -1} + \frac{3}{2}\alpha (t-N)^{\alpha-1} >0 .$$

Since $F'(t) >0 $ , and 
$$\lim_{t \to (N+ 1/2)^-} F(t) = \lim_{t \to (N + 1/2)^+} F(t),$$
so that $F(t)$ reaches minimum when $t=1$ on the domain $t \in (1,\infty)$.

When $N \le -2$, we have 
$$F(t)= [-N + (\frac{1}{2}-t + N)^\alpha]+ [-N + (1-t+N)^\alpha] + [-N-2 + (-1-t - (-N-2)^\alpha)]\frac{3}{2},$$
when $t \in (N,N+1/2),$
and 
$$F(t)= [-N-1 + (\frac{1}{2}-t - (-N-1))^\alpha]+ [-N + (1-t+N)^\alpha] + [-N-2 + (-1-t - (-N-2)^\alpha)]\frac{3}{2},$$
when $t \in (N+1/2,N+1).$
Similar calculation gives $F(t)$ reaches local minimum when $t=-1$ on the domain $t\in (-\infty, -1).$

In the remaining domains, $(0,1/2)$, $(1/2,1)$, $(-1,-1/2)$, $(-1/2,0)$, we may directly calculate the total cost as follows.

When $t\in (0,1/2)$, $$F(t)=[0 + (\frac{1}{2}-t - 0 )^\alpha] + [0 + (1-t-0)^\alpha] + [1 + (1+t -1)^\alpha]\frac{3}{2}.$$

When $t\in (1/2,1)$, $$F(t)=[0 + (t-\frac{1}{2} - 0 )^\alpha] + [0 + (1-t-0)^\alpha] + [1 + (1+t -1)^\alpha]\frac{3}{2}.$$

When $t\in (-1/2,0)$, $$F(t)=[0 + (\frac{1}{2}-t - 0 )^\alpha] + [1 + (1-t-1)^\alpha] + [0 + (1+t -0)^\alpha]\frac{3}{2}.$$

When $t\in (-1,-1/2)$, $$F(t)=[1 + (\frac{1}{2}-t - 1 )^\alpha] + [1 + (1-t-1)^\alpha] + [0 + (1+t -0)^\alpha]\frac{3}{2}.$$

The above four cost functions are concave down, so that potential local minimum values over the domain $t\in (-1,1)$ are $F(-1)$, $F(-1/2)$, $F(0)$, $F(1/2)$, $F(1)$.
Direct calculation gives 
$$F(-1)= 3 + \frac{\sqrt{2}}{2} \approx 3.707, 
\ F(-1/2)= 2 + \frac{\sqrt{2}}{2} + \frac{3\sqrt{2}}{4} \approx 3.767, 
$$
$$F(1)= 3 + \frac{\sqrt{2}}{2} \approx 3.707,
\ F(1/2) = \frac{3}{2} + \frac{\sqrt{2}}{2} + \frac{3\sqrt{2}}{4} \approx 3.267$$
$$F(0)= \frac{5}{2} + \frac{\sqrt{2}}{2} \approx 3.207.$$
Here, $F(-1),\ F(1/2),\ F(1)$ are the costs correspond to the locally cycle-free transport path, which is larger than the cyclic counterpart, $F(0)$.

\end{example}

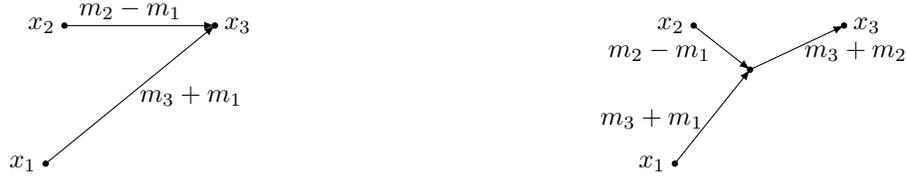
\begin{figure}[h]
\begin{subfigure}{0.45\textwidth}
\centering
\begin{tikzpicture}[>=latex]

\filldraw[black] (0,0) circle (1pt) node[anchor=east]{$x_1$};
\filldraw[black] (0.25,1.837) circle (1pt) node[anchor=east]{$x_2$};
\filldraw[black] (2.25,1.837) circle (1pt) node[anchor=west]{$x_3$};

\filldraw[black] (1.125, 1.8) circle (0pt) node[anchor = south]{$m_2-m_1$};
\filldraw[black] (1.125, 0.9) circle (0pt) node[anchor = west]{$m_3+m_1$};

\draw[->] (0.25,1.837) -- (2.25,1.837);
\draw[->] (0,0) -- (2.25,1.837);

\end{tikzpicture}
\caption{Cycle-free case that corresponds to $T+\frac{1}{2}R.$}
\end{subfigure}
\begin{subfigure}{0.54\textwidth}
\centering
\begin{tikzpicture}[>=latex]

\filldraw[black] (0,0) circle (1pt) node[anchor=east]{$x_1$};
\filldraw[black] (0.25,1.837) circle (1pt) node[anchor=east]{$x_2$};
\filldraw[black] (2.25,1.837) circle (1pt) node[anchor=west]{$x_3$};
\filldraw[black] (1,1.25) circle (1pt) node[anchor=west]{ };

\filldraw[black] (0.6, 1.5) circle (0pt) node[anchor = east]{$m_2-m_1$};
\filldraw[black] (0.5, 0.6) circle (0pt) node[anchor = east]{$m_3+m_1$};
\filldraw[black] (1.6, 1.5) circle (0pt) node[anchor = west]{$m_3+m_2$};

\draw[->] (0.25,1.837) -- (1,1.25);
\draw[->] (0,0) -- (1,1.25);
\draw[->] (1,1.25) -- (2.25,1.837);

\end{tikzpicture}
\caption{Optimizing $T+\frac{1}{2}R$ to a potential Y shaped transport path.}    
\end{subfigure}

\caption{$supp(T+\frac{1}{2}R) \not\subseteq supp(T)$ }
\label{fig: potential acyclic transport path}

\end{figure}

Note that the above calculation is based on assumption or construction that $$supp(T+tR) \subseteq supp(T),$$
which serves as a tool to simplify the overall calculation. 
That being said we may still reach an acyclic result when using $\mathbf{M}_{\alpha,c}$ cost, as indicated by Figure \ref{fig: potential acyclic transport path}, using the labels from the above example.

Proposition \ref{prop: integer and decimal acyclic components} shows how paths with weight equal integer multiple of $c$ interact with each other, 
and Theorem \ref{thm: integer part and decimal part} gives a way to decompose a general transport path into the sum of a transport path with integer multiple of $c$ and a transport path with ``decimal'' multiple of $c$, i.e. weight less than $c$.
In the following, we will show, under the $\mathbf{M}_{\alpha,c}$ cost, transport paths with weight equal to integer multiple of $c$ rarely interact with paths that have weight less than $c$.

\begin{proposition}\label{prop: integer weighted line segments replacement}
Given $\mu^+,\mu^-$ two measures of equal mass, $\alpha \in [0,1]$, $c >0$, and $T= \underline{\underline{\tau}}(M,\theta(x),\xi(x))  \in Path(\mu^-,\mu^+)$. 
Let $v_1,v_2$ be points in $M$ such that 
\begin{equation}\label{eqn: curves between two points}
\Gamma(v_1,v_2):= \{ \gamma(t):[0,1] \to M \ | \ \gamma(0)=v_1, \gamma(1) = v_2, \gamma'(t) \not= 0, \text{ for } t \in (0,1)\}    
\end{equation}
is non-empty. 
Let $\gamma_0 \in \Gamma(v_1,v_2)$ where 
$$\theta_0 = \min \left\{ \left\lfloor \frac{\theta(\gamma_0(t))}{c}\right\rfloor : t \in [0,1] \right\} \ge 1,$$
then 
\begin{equation}\label{eqn: replace curve by edge}
\mathbf{M}_{\alpha,c}(T + \theta_0 c \cdot \llbracket \overline{v_1v_2}\rrbracket - \theta_0 c\cdot \llbracket \gamma_0([0,1])\rrbracket ) \le \mathbf{M}_{\alpha,c}(T),    
\end{equation}
where 
$$\llbracket \overline{v_1v_2}\rrbracket =  \underline{\underline{\tau}}\left(\overline{v_1v_2},1,\frac{v_2 - v_1}{|v_2 - v_1|} \right), \text{ and } 
\llbracket \gamma_0([0,1])\rrbracket =  \underline{\underline{\tau}}\left(\gamma_0([0,1]),1,\frac{\gamma_0'(t)}{|\gamma_0'(t)|} \right).$$
Moreover, equality holds in (\ref{eqn: replace curve by edge}) if and only if $\overline{v_1v_2} = \gamma_0([0,1])$.

Note that $\overline{v_1v_2}$ stands for the direct line segment from $v_1$ to $v_2$, i.e. $\overline{v_1v_2} = \{(1-t)v_1 + tv_2 :  t \in [0,1]\}$.
Moreover, if $\theta_0 = 0$, then equation (\ref{eqn: replace curve by edge}) holds trivially.
\end{proposition}

\begin{proof}
Without loss of generality we may assume $\theta(x)=0$ for $x \in \overline{v_1v_2}$, so that $\overline{v_1v_2} \subseteq M$.
Since 
$$M = M\backslash \left(\overline{v_1v_2} \cup  \gamma_0([0,1]) \right) 
\ \cup\  
\overline{v_1v_2} 
\ \cup \ 
\gamma_0([0,1]),$$
we have the following equation,
\begin{eqnarray*}
\mathbf{M}_{\alpha,c}(T + \theta_0c \cdot \llbracket \overline{v_1v_2}\rrbracket - \theta_0c \cdot \llbracket \gamma_0([0,1])\rrbracket )  
&=& 
\mathbf{M}_{\alpha,c}(\underline{\underline{\tau}}(M\backslash(\overline{v_1v_2} \cup \gamma_0([0,1])),\theta(x),\xi(x)) )
+  \\
&&
\mathbf{M}_{\alpha,c}(\underline{\underline{\tau}}(\overline{v_1v_2},\theta(x)+\theta_0c,\xi(x)) ) +  \\
&&
\mathbf{M}_{\alpha,c}(\underline{\underline{\tau}}(\gamma_0([0,1]),\theta(x) - \theta_0 c,\xi(x)) ). 
\end{eqnarray*}
By definition of $\mathbf{M}_{\alpha,c}$ cost, we have
\begin{eqnarray*}
\mathbf{M}_{\alpha,c}(\underline{\underline{\tau}}(\overline{v_1v_2},\theta(x)+\theta_0c,\xi(x)) ) 
&=&
c^\alpha \int_{\overline{v_1v_2}} \left\lfloor \frac{\theta(x) + \theta_0c}{c}\right\rfloor + \left( \frac{\theta(x) + \theta_0c}{c} -  \left\lfloor \frac{\theta(x) + \theta_0c}{c}\right\rfloor \right)^\alpha  d\mathcal{H}^1  \\
&=&
c^\alpha \int_{\overline{v_1v_2}} \left\lfloor \frac{\theta(x) }{c}\right\rfloor + \theta_0 + \left( \frac{\theta(x)}{c} -  \left\lfloor \frac{\theta(x)}{c}\right\rfloor \right)^\alpha  d\mathcal{H}^1  \\
&=&
c^\alpha \int_{\overline{v_1v_2}}  \theta_0  d\mathcal{H}^1 
+ 
c^\alpha \int_{\overline{v_1v_2}} \left\lfloor \frac{\theta(x) }{c}\right\rfloor + \left( \frac{\theta(x)}{c} -  \left\lfloor \frac{\theta(x)}{c}\right\rfloor \right)^\alpha  d\mathcal{H}^1,
\end{eqnarray*}
and 
\begin{eqnarray*}
\mathbf{M}_{\alpha,c}(\underline{\underline{\tau}}(\gamma_0([0,1]),\theta(x) - \theta_0 c,\xi(x)) )
&=&
c^\alpha \int_{\gamma_0([0,1])} \left\lfloor \frac{\theta(x) - \theta_0c}{c}\right\rfloor + \left( \frac{\theta(x) - \theta_0c}{c} -  \left\lfloor \frac{\theta(x) - \theta_0c}{c}\right\rfloor \right)^\alpha  d\mathcal{H}^1  \\
&=&
c^\alpha \int_{\gamma_0([0,1])} \left\lfloor \frac{\theta(x) }{c}\right\rfloor - \theta_0 + \left( \frac{\theta(x)}{c} -  \left\lfloor \frac{\theta(x)}{c}\right\rfloor \right)^\alpha  d\mathcal{H}^1  \\
&=&
-c^\alpha \int_{\gamma_0([0,1])}  \theta_0  d\mathcal{H}^1 
+ 
c^\alpha \int_{\gamma_0([0,1])} \left\lfloor \frac{\theta(x) }{c}\right\rfloor + \left( \frac{\theta(x)}{c} -  \left\lfloor \frac{\theta(x)}{c}\right\rfloor \right)^\alpha  \!\!d\mathcal{H}^1.
\end{eqnarray*}
Hence, 
\begin{eqnarray*}
\mathbf{M}_{\alpha,c}(T + \theta_0 \cdot \llbracket \overline{v_1v_2}\rrbracket - \theta_0 \cdot \llbracket \gamma_0([0,1])\rrbracket )  - \mathbf{M}_{\alpha,c}(T) 
&=& 
c^\alpha \int_{\overline{v_1v_2}}  \theta_0  d\mathcal{H}^1  -
c^\alpha \int_{\gamma_0([0,1])}  \theta_0  d\mathcal{H}^1  \\ 
&=&
c^\alpha \theta_0 \cdot \left( \mathcal{H}^1 (\overline{v_1v_2}) - \mathcal{H}^1(\gamma_0([0,1])) \right).
\end{eqnarray*}

Since $\gamma_0 \in \Gamma(v_1,v_2)$, we have $\gamma_0(0) = v_1, \gamma_0(1) = v_2$. 
Moreover, the $1$ dimensional Hausdorff measure of a curve equals the length of the curve, which implies $$\mathcal{H}^1 (\overline{v_1v_2}) - \mathcal{H}^1(\gamma_0([0,1])) \le 0,$$
and the equality takes place if and only if $\overline{v_1v_2} = \gamma_0([0,1])$.
Since we assume $\theta_0 \ge 1$, we get inequality (\ref{eqn: replace curve by edge}).
\end{proof}

\begin{corollary}\label{cor: optimal when min density larger than c}
Given $\mu^+,\mu^-$ two measures of equal mass, $\alpha \in [0,1]$, $c >0$, and $T= \underline{\underline{\tau}}(M,\theta(x),\xi(x))  \in Path(\mu^-,\mu^+)$ be an optimal transport path under $\mathbf{M}_{\alpha,c}$ cost.
Let $v_1,v_2$ be points in $M$ such that 
$\Gamma(v_1,v_2)$ (defined in (\ref{eqn: curves between two points})) is non-empty.
Suppose $\gamma_0 \in \Gamma(v_1,v_2)$ such that 
$$\theta_0 = \min \left\{  \frac{\theta(\gamma_0(t))}{c} : t \in [0,1] \right\} \in \mathbb{Z}^+ ,$$
then we have $\gamma_0([0,1]) = \overline{v_1v_2}$.
\end{corollary}

\begin{proof}
Proof by contradiction and suppose $\gamma_0([0,1]) \not= \overline{v_1v_2}$, then by Proposition \ref{prop: integer weighted line segments replacement} we have 
$$\mathbf{M}_{\alpha,c}(T + \theta_0 c \cdot \llbracket \overline{v_1v_2}\rrbracket - \theta_0 c\cdot \llbracket \gamma_0([0,1])\rrbracket ) 
<
\mathbf{M}_{\alpha,c}(T).$$
Since $$\partial  \left( T + \theta_0 c \cdot \llbracket \overline{v_1v_2}\rrbracket - \theta_0 c\cdot \llbracket \gamma_0([0,1])\rrbracket \right) = \partial T + \partial \left( \theta_0 c \cdot \llbracket \overline{v_1v_2}\rrbracket - \theta_0 c\cdot \llbracket \gamma_0([0,1])\rrbracket \right) = \partial T,$$
we have $T + \theta_0 c \cdot \llbracket \overline{v_1v_2}\rrbracket - \theta_0 c\cdot \llbracket \gamma_0([0,1])\rrbracket \in Path(\mu^-,\mu^+)$, which contradicts the assumption that $T$ is an optimal transport path under the cost $\mathbf{M}_{\alpha,c}$, i.e. $\mathbf{M}_{\alpha,c}(T) \le \mathbf{M}_{\alpha,c}(\tilde{T})$ for all $\tilde{T} \in  Path(\mu^-,\mu^+)$.
\end{proof}

\begin{proposition}\label{prop: integer components separation}
Given $\mu^+,\mu^-$ two measures of equal mass, $\alpha \in [0,1]$, $c >0$, and $T= \underline{\underline{\tau}}(M,\theta(x),\xi(x))  \in Path(\mu^-,\mu^+)$. 
Let $v_1,v_2$ be points in $M$ such that 
$\Gamma(v_1,v_2)$ (defined in (\ref{eqn: curves between two points})) is non-empty.
Suppose $\gamma_0 \in \Gamma(v_1,v_2)$ such that 
$$\theta_0 = \min \left\{ \frac{\theta(\gamma_0(t))}{c} : t \in [0,1] \right\} \ge 1.$$
For arbitrary $\theta_1 \in \mathbb{Z}^{\ge 0}$, with $\theta_1 \le \theta_0$, 
we have 
$$\mathbf{M}_{\alpha,c}(\underline{\underline{\tau}}(\gamma_0([0,1]),\theta(x) - \theta_1 c,\xi(x)) ) 
+ 
\mathbf{M}_{\alpha,c}(\underline{\underline{\tau}}(\gamma_0([0,1]),\theta_1 c,\xi(x)) ) 
=
\mathbf{M}_{\alpha,c}(\underline{\underline{\tau}}(\gamma_0([0,1]),\theta(x),\xi(x)) ).$$
This gives that under the cost function $\mathbf{M}_{\alpha,c}$, we may identify the curve $\gamma_0([0,1])$ with weight $\theta(x)$ as two curves $\gamma_0([0,1])$ with weight $\theta(x) - \theta_1 c$ and $\theta_1 c$ respectively.

\end{proposition}

\begin{proof}
The proof is just a direct calculation of $\mathbf{M}_{\alpha,c}$ cost, which is the same as the calculation in the proof of Proposition \ref{prop: integer weighted line segments replacement}.
\begin{eqnarray*}
\mathbf{M}_{\alpha,c}(\underline{\underline{\tau}}(\gamma_0([0,1]),\theta(x) - \theta_1 c,\xi(x)) )
&=&
c^\alpha \int_{\gamma_0([0,1])} \left\lfloor \frac{\theta(x) - \theta_1 c}{c}\right\rfloor + \left( \frac{\theta(x) - \theta_1 c}{c} -  \left\lfloor \frac{\theta(x) - \theta_1 c}{c}\right\rfloor \right)^\alpha  d\mathcal{H}^1  \\
&=&
c^\alpha \int_{\gamma_0([0,1])} \left\lfloor \frac{\theta(x) }{c}\right\rfloor - \theta_1 + \left( \frac{\theta(x)}{c} -  \left\lfloor \frac{\theta(x)}{c}\right\rfloor \right)^\alpha  d\mathcal{H}^1  \\
&=&
-c^\alpha \int_{\gamma_0([0,1])}  \theta_1  d\mathcal{H}^1 
+ 
c^\alpha \int_{\gamma_0([0,1])} \left\lfloor \frac{\theta(x) }{c}\right\rfloor + \left( \frac{\theta(x)}{c} -  \left\lfloor \frac{\theta(x)}{c}\right\rfloor \right)^\alpha  d\mathcal{H}^1 \\
&=&
-\mathbf{M}_{\alpha,c}(\underline{\underline{\tau}}(\gamma_0([0,1]),\theta_1 c,\xi(x)) ) 
+
\mathbf{M}_{\alpha,c}(\underline{\underline{\tau}}(\gamma_0([0,1]),\theta(x),\xi(x)) ).
\end{eqnarray*}
\end{proof} 

Using Corollary \ref{cor: optimal when min density larger than c} and Proposition \ref{prop: integer components separation}, we may characterize optimal transport paths under $\mathbf{M}_{\alpha,c}$ cost in some simple cases. 

\begin{example}
Let $\mu^- = 2.5\delta_{x_1} + 0.5\delta_{x_2}$, $\mu^+ = 3 \delta_y$, where $x_1 = (-1,3)$, $x_2 = (1,3)$, $y=(0,0)$, and $c=1$.

\begin{figure}[h]
\begin{subfigure}{0.32\textwidth}
\centering
\begin{tikzpicture}[>=latex]
\filldraw[black] (-1,3) circle (1pt) node[anchor=east]{$x_1$};
\filldraw[black] (1,3) circle (1pt) node[anchor=west]{$x_2$};
\filldraw[black] (0,0) circle (1pt) node[anchor=east]{$y$};
\filldraw[black] (0,1.8) circle (0pt) node[anchor = east]{ };

\filldraw[black] (-0.6, 2.5) circle (0pt) node[anchor = north]{ };
\filldraw[black] (0.65, 2.5) circle (0pt) node[anchor = north]{ };
\filldraw[black] (0, 1) circle (0pt) node[anchor = east]{ };


\end{tikzpicture}
\caption{}
\end{subfigure}
\begin{subfigure}{0.32\textwidth}
\centering
\begin{tikzpicture}[>=latex]
\filldraw[black] (-1,3) circle (1pt) node[anchor=east]{$x_1$};
\filldraw[black] (1,3) circle (1pt) node[anchor=west]{$x_2$};
\filldraw[black] (0,0) circle (1pt) node[anchor=east]{$y$};
\filldraw[black] (0,1.8) circle (0pt) node[anchor = east]{ };

\filldraw[black] (-0.5, 1.5) circle (0pt) node[anchor = east]{$2$};
\filldraw[black] (0.65, 2.5) circle (0pt) node[anchor = north]{ };
\filldraw[black] (0, 1) circle (0pt) node[anchor = east]{ };

\draw[thick,->] (-1,3) -- (0,0);

\end{tikzpicture}
\caption{}
\end{subfigure}
\begin{subfigure}{0.32\textwidth}
\centering
\begin{tikzpicture}[>=latex]
\filldraw[black] (-1,3) circle (1pt) node[anchor=east]{$x_1$};
\filldraw[black] (1,3) circle (1pt) node[anchor=west]{$x_2$};
\filldraw[black] (0,0) circle (1pt) node[anchor=east]{$y$};
\filldraw[black] (0,2) circle (1pt) node[anchor = east]{ };
\filldraw[black] (0,1.8) circle (0pt) node[anchor = east]{ };

\filldraw[black] (-0.5, 1.5) circle (0pt) node[anchor = east]{$2$};
\filldraw[black] (-0.6, 2.6) circle (0pt) node[anchor = west]{$0.5$};
\filldraw[black] (0.5, 2.5) circle (0pt) node[anchor = west]{$0.5$};
\filldraw[black] (0, 1) circle (0pt) node[anchor = west]{$1$};

\draw[thick,->] (-1,3) -- (0,0);
\draw[->] (-1,3) -- (0,2);
\draw[->] (1,3) -- (0,2);
\draw[->] (0,2) -- (0,0);

\end{tikzpicture}
\caption{}
\end{subfigure}

\caption{ }
    
\label{fig: simple optimal cases }
\end{figure}

Since $\mu^+$ is supported on only $1$ point, all the weight from $x_1$ and $x_2$ will be transported to $y$ via some curve or $1-$current. 
By Proposition \ref{prop: integer components separation},
we may identify there are two curves starting from $x_1$, one curve has weight $2$ and the other has $0.5$.
By Corollary \ref{cor: optimal when min density larger than c}, we may assume the curve with density $2$ is a straight line segment $\overline{x_1 y}$, which is the second graph above.

Since the remaining total weight from $x_1$ and $x_2$ is less or equal to $1$, results from either Remark \ref{rem: 2points to 1 point calculation} or \cite{xia2003} gives the ``Y'' shaped paths indicated above in the last graph.

\end{example}

\end{document}